\documentclass[lang=american, preliminary]{ems-icm}

\usepackage{framed,bm}

\numberwithin{equation}{section}

\begin{document}

\firstpage{1}
\doinumber{204}
\volume{1}
\copyrightyear{2022}

\title{100 Years of the (Critical) Ising Model on the Hypercubic Lattice}
\titlemark{100 Years of the Ising Model}

\emsauthor{1}{Hugo Duminil-Copin}{Hugo Duminil-Copin}


\emsaffil{1}{Unige, 
7-9 rue du Conseil-G\'en\'eral,
1205 Gen\`eve, Switzerland~\email{hugo.duminil@unige.ch},
IHES, 35 route de Chartres, 91440 Bures-Sur-Yvette, France \email{duminil@ihes.fr}.}

\dedication{Dedicated to the memory of colleagues and friends Dmitry Ioffe and Vladas Sidoravicius}

\classification[60K35]{82B20}

\keywords{Ising model, percolation theory, phase transition, statistical mechanics}

\begin{abstract}
We take the occasion of this article to review one hundred years of the physical and mathematical study of the Ising model. The model, introduced by Lenz in 1920, has been at the cornerstone of many major revolutions in statistical mechanics. We wish, through its history, to outline some of these amazing developments. We restrict our attention to the ferromagnetic nearest-neighbour model on the hypercubic lattice, and essentially focus on what happens at or near the so-called critical point.
\end{abstract}

\maketitle


\section{Short motivation} \label{sec:1}

How to provide an introduction to (part of) statistical physics aimed at a wide audience of mathematicians? The question is not easy, especially since the domain is positioned halfway between theoretical physics and mathematics, and that contrarily to (some) other fields of mathematics, it is hard to identify a theory that would embrace most of statistical physics. A partial answer may be to follow the standard approach of teaching by examples, and to pick what is maybe the most classical model of statistical physics, namely the Ising model. Through its history, one may trace many of the revolutions, both on the theoretical physics and mathematical sides, that statistical physics underwent in the last century. 

We therefore chose to streamline this history, from the emergence of the model to explain experimental results, to its modern applications in 
mathematics, physics, and beyond. Obviously, the story will be tainted by the expertise of the author since thousands of papers have been published on the subject, which ranges over many subfields of mathematics and theoretical physics. A subjective selection of papers has therefore been made, and the attention has been restricted to the model  on the hypercubic lattice, at equilibrium, (for most of the review) at criticality, and always with nearest neighbour ferromagnetic interactions (we are well aware that non-critical and dynamical aspects, as well as long-range, random, or antiferromagnetic interactions also are of prime importance).

We tried to respect the timing of the appearance of the different notions pertaining to the model, and avoided as much as possible some tempting anachronisms. As a result, certain readers may be surprised by some statements, knowing that simpler and more natural versions exist nowadays. Also, the large number of breakthroughs in the Ising model's history -- Peierls' argument, Onsager's solution on the square lattice and the exact integrability results that followed, Kadanoff's scaling and universality hypotheses, the Lee-Yang theorem, correlation inequalities, the Fortuin-Kasteleyn representation, reflection positivity, Aizenman's treatment of the random current representation and use of differential inequalities, conformal field theory, rigorous renormalisation group, Chelkak-Smirnov's conformal invariance, 3D conformal bootstrap, to cite but a few -- forced us to be very quick on some of these developments. References are added for the avid reader. We also refer to \cite{Bru67,Nis} for historical introductions, and \cite{FriVel17} for a book on statistical physics including a study of the Ising model.

\section{The first twenty years: a laborious start} \label{sec:2}

\subsection{Ising model's prehistory} \label{sec:2.1}

In 1895, the French physicist Pierre Curie \cite{Cur85} noticed that magnets lose their magnetic attraction when they are heated above a certain critical temperature, now called the {\em Curie temperature} (the phenomenon seemed to have been  discovered before by the French physicist Pouillet in 1832, see \cite{Jon} for precise references and a discussion). While the Curie temperature varies from slightly over 100 degrees Celsius for certain alloys, to 769.85 degrees Celsius for magnets made of iron, the underlying phenomenon is always the same: at a certain temperature, a magnet ceases to be able to keep a {\em spontaneous magnetisation} and exhibits magnetisation {\em only} when an external field is applied to it. This phenomenon is called a {\em phase transition} between a paramagnetic phase above the Curie temperature, and a ferromagnetic phase below it.

Curie also identified a law, now called {\em Curie's law}, relating the magnetic susceptibility of the system to the temperature and the external magnetic field applied to the magnet. He noticed the similarity between the ferromagnetic and paramagnetic phases of a magnet in terms of the temperature and the external magnetic field applied to it, and respectively the liquid and gas phases of a fluid in terms of the pressure and temperature. Pierre Weiss \cite{Wei07} tried to produce an efficient physical explanation of this phenomenon by introducing an assumption, referred to today as the {\em mean-field approximation}. This mean-field model, called the {\em Curie-Weiss model}, gave rise to an interesting yet not fully accurate description of the phase transition. 

The German physicist Wilhelm Lenz got interested in Curie's law. Lenz agreed with one of Weiss' suggestions that magnets are made of elementary pieces that behave themselves as small magnets. Yet, he was at the same time in line with his contemporary physicists, thinking that one of Weiss' assumption, namely that elementary magnets can rotate freely within a solid, was wrong. Taking this into account, he challenged the rotational freeness. Observing that a crystal selects certain directions corresponding to its symmetries, he made
 the assumption that elementary magnets also behave in this way. By analogy, he then suggested that a crystal-like mechanism for magnets should favour that neighbouring elementary magnets are aligned, therefore corresponding to either pointing in the same or opposite directions. At the end, the reasoning of Lenz led to the assumption that elementary magnets were taking only two possible directions that are opposite of each other. He formalized this reasoning in \cite{Len20}.

At this stage, Lenz did not propose an explicit form for the interaction between elementary magnets. Also, the paper approximately explained the typical behaviour of a paramagnet having respectively zero magnetisation when no magnetic field is applied, and a magnetisation when such a magnetic field is applied, but Lenz made no mention of what will later be referred to as the ferromagnetic behaviour.

Ernst Ising was a German physicist born in 1900, who was a PhD student of Lenz in Hamburg. He graduated in 1924 and published a paper \cite{Isi25} on Lenz's model in 1925. 
So, what did Ising actually achieve in his famous paper from 1925? 

First of all, he went one step further than Lenz by specifying the interaction between elementary magnets. He first made the assumption that interactions "{\em decay rapidly with distance, so that we, in general, to a first approximation, only have to take the influence on neighbouring elements into account}". He also assumed that "{\em of all the possible positions that the neighbouring atoms can assume in relation to each other, the one that requires the minimum energy is when they are both acting in the same direction}". These two assumptions led to the mathematical model that we will define formally in the next section. In order to treat this model, Ising made a further assumption: he assumed that the elementary magnets are positioned on a linear chain. 

From all of this, Ising could deduce Curie's law in the paramagnetic phase. While this was a source for optimism, the latter was severely challenged by the observation that the magnetisation was tending to 0 as the magnetic field vanishes, irrespective of the temperature. In other words, no explanation for ferromagnetism was in sight. Even worse, despite a few attempts at generalising the model (Ising considered non-nearest neighbour interactions, more possible directions for the elementary magnets, and a hybrid three dimensional model that would correspond to a limit in which only pairs of neighbouring elementary magnets in one direction truly interact), the ferromagnetism did not seem to be explainable. This led Ising to conjecture that the model was not a good explanation for ferromagnetism (even when considering higher dimensional base graph for the spins), a thought that he gathered in a letter to American historian Stephen Brush years later:
"{\em I discussed the result of my paper widely with Professor Lenz and with Dr. Wolfgang Pauli, who at that time was teaching in Hamburg. There was some disappointment that the linear model did not show the expected ferromagnetic properties}."

After his PhD, Ising left academia to become a teacher in Germany before being forced to step down due to his Jewish origins. He fled Nazi Germany and emigrated to the United States, where he became a Professor in Physics at Bradley University. He never published after his first original paper, and only later became aware of how famous the model had grown into.

\subsection{Formal definition} \label{sec:2.2}

Let us turn to the formal definition of the model for our magnet. Consider a finite non-oriented subgraph $G=(V,E)$ of the hypercubic lattice $\mathbb Z^d$ with vertex-set $V$ corresponding to the position of its elementary magnet constituents, and edge-set $E$ modeling the links between neighbouring ones. An edge $e\in E$ is often written $e=\{x,y\}$, where $x$ and $y$ are its endpoints. The elementary magnet at $x\in V$ will be a quantity $\sigma_x\in \{-1,+1\}$, where $-1$ and $+1$ correspond to the two opposite directions that it may take. The value $\sigma_x$ is called the {\em spin at $x$}, and the collection $(\sigma_x:x\in V)\in\{-1,1\}^V$ of all spins at vertices in $V$ is called the {\em spin configuration}, and should be understood as the state of our magnet.

The {\em energy} -- or {\em Hamiltonian} -- of a configuration $\sigma$ on $G$ is given by 
\begin{equation}
H_{G,h}(\sigma):=-\sum_{\{x,y\}\in E}\sigma_x\sigma_y-h\sum_{x\in V}\sigma_x,
\end{equation}
where $h\in \mathbb R$ is called the {\em magnetic field}. Sometimes, one may want to generalise the model to accommodate non-nearest neighbour and non-ferromagnetic interactions by setting
\begin{equation}\label{eq:def1}
H_{G,h,(J_{x,y})}(\sigma):=-\sum_{x,y\in V}J_{x,y}\sigma_x\sigma_y-h\sum_{x\in V}\sigma_x
\end{equation}
where the $(J_{x,y}:x,y\in V)$ are called the {\em coupling constants} of the model. Except when otherwise stated, we focus here on the Hamiltonian $H_{G,h}$ corresponding to what is called the {\em nearest neighbour ferromagnetic} (n.n.f.) Ising model on $G$. 

Following Boltzmann, one considers the {\em (grand) partition function} of the Ising model on $G$ at inverse-temperature $\beta$ and magnetic field $h$ defined by
\begin{equation}\label{eq:def2}
Z(G,\beta,h):=\sum_{\sigma\in\{-1,1\}^V}\exp[-\beta H_{G,h}(\sigma)].
\end{equation}
The quantity $\beta$ is interpreted as the inverse of the temperature, as the latter corresponds to the thermal excitation of elementary magnets, for which it is natural to predict that the larger their excitation, the less relevant their interaction. 

Physicists then consider the linear form defined for every function $X:\{-1,1\}^V\rightarrow \mathbb R$ by the formula
\begin{equation}
\langle X\rangle_{G,\beta,h}:=\frac1{Z(G,\beta,h)}\sum_{\sigma\in \{-1,1\}^V}X(\sigma)\exp[-\beta H_{G,\beta,h}(\sigma)].
\end{equation}

At this stage, we do not consider $\langle\cdot\rangle_{G,\beta,h}$ itself, and instead focus on a thermodynamical quantity of the system called the free energy.
 Consider a box $\Lambda_n:=[-n,n]^d\cap\mathbb Z^d$ and define the {\em free energy} of the $d$-dimensional Ising model by the formula 
  \begin{equation}\label{eq:free energy}
  f(\beta,h):=-\frac1\beta\lim_{n\rightarrow \infty} \frac1{|\Lambda_n|}\ln Z(\Lambda_n,\beta,h)
  \end{equation}
(the existence of the limit is justified by a sub-additivity argument left to the reader). 

Originally, Lenz and Ising were interested in a quantity 
\begin{equation}
m(\beta,h):=-\tfrac{\partial}{\partial h}f(\beta,h),
\end{equation}
which is interpreted as the {\em magnetisation} of the system in the presence of a magnetic field of strength $h$. One may then define the {\em spontaneous magnetisation}, which  corresponds to the remaining magnetisation when removing the magnet from the ambient magnetic field, 
\begin{equation}
m^*(\beta):=\lim_{h\searrow 0}m(\beta,h)
\end{equation} 
(to justify the limit, one may prove that $m(\beta,h)$ decreases as $h$ decreases). The cases $m^*(\beta)=0$ and $m^*(\beta)>0$ are respectively called the paramagnetic and ferromagnetic cases as they correspond to the cases where the magnet respectively loses or keeps its magnetisation even without external magnetic field. 

\subsection{What does the Ising model truly model?
} \label{sec:2.3}

The Ising model did not develop quickly after its introduction. The original paper was cited very sporadically in the ten years that followed. In fact, Ising himself was aware of one citation to his paper only, and this lack of interest was one of the reasons that pushed him to abandon academia. 

There are several explanations why the paper received little attention. The first one is that the negative result of the paper, stating that the model did not explain ferromagnetism, was a pretty disappointing one. The second is a timing problem. A few years after Ising's paper, Heisenberg introduced another model of ferromagnetism \cite{Hei28} based on quantum mechanics, in which the ``classical'' spins of the Ising model are replaced by the quantum spins of electrons. In other words, Heisenberg's model tries to explain ferromagnetism via the interaction of the spin angular momentum of the electrons in the atoms, while the Ising model was relying on their magnetic moments. In a certain sense, the Ising model was a semi-classical version of Heisenberg model, and as such was violating the latest developments of quantum mechanics. The discrepancy between the great predictive successes of the Heisenberg model, and the impossibility to reconcile the Ising model with the recent advancements in modern physics almost entirely disqualified the model as a good description of ferromagnetic materials.

At this point, one may wonder why this model, initially introduced in theoretical physics to explain ferromagnetism but seemingly failed to do so, did not simply fall into darkness after this rocky start. An element of answer can be found in the developments of other fields of physics, which we now review.

In 1919, the Russian-German chemical physicist Gustav Tamman presented an interesting experiment in which atoms in alloys of copper and gold tend to be surrounded by atoms of the other kind (to picture this, think of a chessboard colouring of the square lattice). In Tamman's experiment, the thermal agitation has a direct impact on how much the atoms tend to be in the right places. In 1935, Bragg and Williams \cite{BraWil34} explained this phenomenon by a statistical mechanics's argument involving the energy cost of having an atom in the wrong place. Hans Bethe simplified the model by assuming that only nearest atoms interact. 

In 1936, Ralph Fowler and his team in Cambridge introduced another theoretical model to understand the adsorption of metal vapour on a glass. Fowler more generally identified a class of experiments exhibiting similar behaviours, that he named {\em cooperative phenomena}.  

The German theoretical physicists Rudolf Peierls
later noticed the similarity between Bethe's approximation of the Bragg-Williams model, Fowler's theory of adsorption, and the Ising model. While the original physical problems are different, the mathematical treatment is in fact similar. In retrospect, Peierls was perhaps the first person to identify that the Ising model could treat a number of different phenomena, even though the model was a coarse caricature for each one of the phenomena in question.

This observation was maybe what kept the Ising model alive for some years, but it is mathematics that truly changed the nature of the model and made it what it is today. We now turn to the first mathematical breakthrough in the model.

\subsection{Peierls' argument} \label{sec:2.4}

While Peierls agreed with the majority of the physics community that the Ising model was not a good model for ferromagnetism, he certainly recognised that the model was of mathematical interest. Furthermore, he totally disagreed with the naive generalisation, based on the few attempts of Ising, of the absence of a ferromagnetic phase to higher dimensional lattices. This led him to reconsider the problem of the Ising model in two and three dimensions. As a result, he produced what is probably one of the most important papers in the early Ising history \cite{Pei36}, in which he developed a technique which is now widely known in statistical physics as {\em Peierls' argument}.

Roughly speaking, the argument runs as follows. When considering a configuration $\sigma$ of the Ising model on $\mathbb Z^2$, or a finite subgraph of it, one may 
 associate a subset $E(\sigma)$ composed of the edges $\{x,y\}$ of the graph  with $\sigma_x\ne \sigma_y$. In a planar context, one may draw these edges $e\in E(\sigma)$ by considering the dual edges $(e^*:e\in E(\sigma))$ on the dual graph\footnote{The dual graph $G^*=(V^*,E^*)$ of a planar graph $G=(V,E)$ is the planar graph with vertex-set given by the faces of $G$ (including the exterior one) and edge-set $E^*$ given by unordered pairs $\{u,v\}$, where $u$ and $v$ are two faces that are bordered by the same edge. When this edge is $e$, we denote the dual edge $\{u,v\}$ by $e^*$. The map $e\mapsto e^*$ is therefore a bijection between $E$ and $E^*$. On the square lattice, the dual graph is nothing but the translate by $(\tfrac12,\tfrac12)$ of the square lattice itself.}; see Figure~\ref{fig:1}. These dual edges and their endpoints form an even subgraph of the dual graph (call $\mathrm{Even}(G^*)$ the set of such even subgraphs) which can be interpreted as a collection of loops on the dual graph. The representation of configurations $\sigma$ in terms of even subgraphs is called the {\em low-temperature expansion}.
  Using the mapping between $\sigma$ and $E(\sigma)$, one may rewrite the partition function as
  \begin{equation}\label{eq:KW1}
  Z(G,\beta,0)=\sum_{\sigma\in\{-1,1\}^{V}} e^{-\beta H_{G,h}(\sigma)}=e^{\beta |E(G)|}\sum_{F\in \mathrm{Even}(G^*)}e^{-2\beta|F|}.
  \end{equation}
 This formula immediately highlights the fact that $\beta$ large renders configurations $F\in \mathrm{Even}(G^*)$ with large loops unlikely. Building on this observation, Peierls was able to obtain that $m^*(\beta)>0$ for large values of $\beta$, see Frame~1 for more details.
 
   The idea to introduce a model of ``domain walls'' separating the different phases (here pure $+1$ and pure $-1$) from each other is not restricted to the Ising model: it has been very fruitful to prove the existence of phase transitions, and Peierls' argument is now one of the most famous and robust arguments in statistical physics.
     \begin{framed}
  \centerline{\textbf{Frame 1: a quick version of Peierls' argument}}
We do not consider the magnetization $m^*(\beta)$ but rather the correlation $\langle\sigma_0\sigma_{\mathbf{g}}\rangle_{\Lambda_n^+,\beta,0}$, where $\Lambda_n^+$ is the graph $\Lambda_n$ plus a vertex $\mathbf{g}$, sometimes referred to as Griffiths' ``ghost'' vertex, connected to all the vertices on the boundary of $\Lambda_n$; see Figure~\ref{fig:1} on the left. The limit as $n$ tends to infinity can be shown to be $m^*(\beta)$, so it is sufficient to prove that the quantity is bounded away from 0 uniformly in $n$.

If one denotes by $\mathbf C=\mathbf C(\sigma)$ the connected component of $0$ in $\mathbb R^2\setminus\{e^*:e\in E(\sigma)\}$, one may decompose the magnetisation depending on the value of $\mathbf C$ to get
 \begin{equation}\label{eq:peierls1}
\langle\sigma_0\sigma_{\mathbf{g}}\rangle_{\Lambda_n^+,\beta,0}=1-2\sum_{\mathbf g\notin C\in \mathrm{Even}((\Lambda_n^+)^*)}\langle \mathbb I(\mathbf C=C)\rangle_{\Lambda_n^+,\beta,0}.
  \end{equation}
Now things become interesting. For every $C\notni \mathbf g$, consider the configuration $\mathrm{Flip}_C(\sigma)$ obtained from $\sigma$ by flipping the values of the spins inside $C$. This effectively corresponds to removing the set $\partial_e C$ of edges in $E(\sigma)$ with {\em exactly one} endpoint in $C$. Taking into account the cost of this operation leads to 
 \begin{align*}
\langle \mathbb I(\mathbf C=C)\rangle_{\Lambda_n^+,\beta,0}\le e^{-2\beta|\partial_e C|}
 \end{align*}
for every $C\notni \mathbf g$. At this stage, the fact that $\partial_eC$ is a loop and that there are at most $(k+1)4^k$ possible loops of length $k$ surrounding the origin  gives
 \begin{equation}\label{eq:peierls1}
\langle\sigma_0\sigma_{\mathbf{g}}\rangle_{\Lambda_n^+,\beta,0}\ge1-2\sum_{k\ge1}k4^ke^{-2\beta k}>1-\frac{8e^{-2\beta}}{(1-4e^{-2\beta})^2}.
  \end{equation}
  \end{framed}

\section{Onsager's 1944 revolution and the integrability of the Ising model} \label{sec:3}

  \subsection{Kramers-Wannier treatment of the Ising model and duality}\label{sec:3.1}
  
  While Peierls' result is certainly one of the first key rigorous steps in the understanding of the Ising model, the work \cite{KraWan41} of Hans Kramers and Gregory Wannier in 1941 propelled the Ising model in another dimension in terms of mathematical interest. Indeed, the two physicists agreed that the Ising model was not necessarily an accurate description of ferromagnetism, but they were precursors in strongly believing that having mathematical models that can be rigorously analysed was of crucial interest for the understanding of physical phenomena, even if only approximate.
  
Kramers and Wannier's goal was to understand what happens for the Ising model at {\em arbitrary} inverse-temperature. Peierls' argument shows that the model behaves like a ferromagnet when $\beta$ is large. A fairly simple argument, see Frame~4, shows that it behaves like a paramagnet when $\beta$ is small. It is therefore tempting to think that there is an intermediate inverse-temperature, playing the theoretical role of the inverse of Curie's temperature, that separates a paramagnetic phase from a ferromagnetic phase, i.e.~a {\em critical inverse-temperature} $\beta_c$ defined by the formula
\begin{equation}
\beta_c=\beta_c(\mathbb Z^d):=\inf\{\beta:m^*(\beta)>0\}.
\end{equation}
Of course, the notion of critical inverse-temperature immediately leads to the following question:  can one compute the value of the critical point $\beta_c$?

  \begin{figure}[t]
  \begin{center} \includegraphics[width=0.90\textwidth]{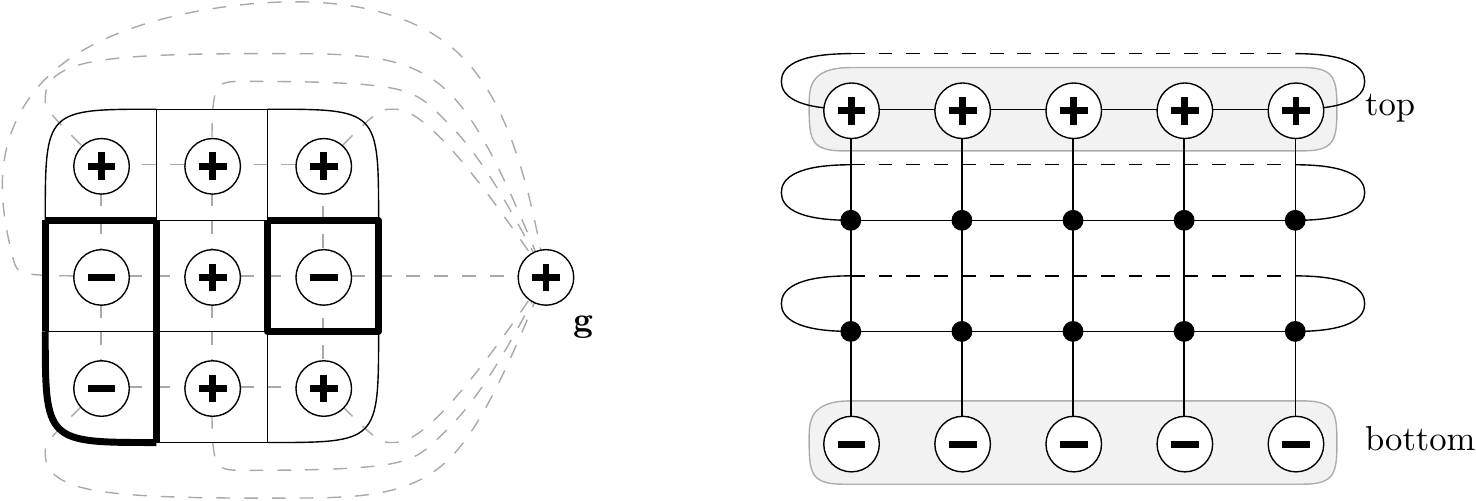} 
\caption{\textbf{On the left.} A picture of the low-temperature expansion on $\Lambda_1^+$. The set $\Lambda_1^+$ is depicted in dashed gray, and $(\Lambda_1^+)^*$ in plain black. The edges in the dual configuration are depicted in bold. \textbf{On the right.}  The set $\mathbb S(5,3)$ with the bottom and top sets depicted. In this case $\underline\tau$ and $\overline\tau$ are respectively constant equal to $-1$ and $+1$.}\label{fig:1}\end{center}\end{figure}

  The work \cite{KraWan41} represented an important historical step towards this computation. It involved a number of ideas that deeply influenced the way mathematicians and physicists approach the Ising model. The first key observation is that Kramers and Wannier did not work with the Ising model in the presence of a magnetic field (in other words, they set $h$ to be 0). Instead, they proposed to look at the {\em specific heat} defined by
  \begin{equation}
  C(\beta):=-\beta^2\tfrac{\partial^2}{(\partial \beta)^2} (\beta f)(\beta,0).
  \end{equation}
  Kramers and Wannier argued that the critical point of the model on $\mathbb Z^2$ should correspond to a value of $\beta$ at which $C(\beta)$ blows up. The next step is maybe the most interesting one. By {\em assuming} that there exists a unique point at which $C(\beta)$ blows up, they were able to predict the value of $\beta_c$. The reason behind this prediction is the following {\em duality relation} for the free energy on $\mathbb Z^2$:
  \begin{equation}
\beta f(\beta,0)=\beta^*f(\beta^*,0)-2\beta+\ln 2+2\ln\cosh(\beta^*),
\end{equation}
where $\beta$ and $\beta^*$ are related via the formula $\tanh(\beta^*)=e^{-2\beta}$. The uniqueness implies that $\beta_c$ must be the self-dual point satisfying $\beta^*=\beta$, i.e. $\beta_c$ must be equal to $\tfrac12\ln(1+\sqrt 2)$. Of course, this reasoning is not a formal proof as it is a priori non-obvious that the singular point is unique.

  The proof of Kramers and Wannier of the duality relation is also of great interest. Originally, they used so-called {\em transfer matrices} to do it; see Frame~2 for details. While they did not invent those matrices (they already appeared in the work of Montroll \cite{Mon41}), they probably made the first important use of them. 
    Today, the derivation of this relation is straightforward and does not rely on transfer matrices. 
 It involves relating the partition functions $Z(G,\beta,0)$ and $Z(G^*,\beta^*,0)$ using, for the first one, the expression given by the low-temperature expansion \eqref{eq:KW1}, and for the second, an alternative representation called the {\em high-temperature expansion}, obtained by van der Waerden \cite{Wae41} and briefly described in Frame~4. When observing that the dual of a box in the square lattice is (except on the boundary) a box of the square lattice, one obtains the identity above by considering larger and larger boxes.

    \begin{framed}
  \centerline{\textbf{Frame 2: transfer matrices of the Ising model}
  }
To lighten the presentation we restrict our attention to the case $h=0$. Consider the slices $\mathbb S(N,M):=(\mathbb Z/N\mathbb Z)^{d-1}\times\llbracket0,M\rrbracket$ with no edges between the vertices of the bottom $(\mathbb Z/N\mathbb Z)^{d-1}\times\{0\}$ (we call $(\mathbb Z/N\mathbb Z)^{d-1}\times\{M\}$ the top of the slice); see Figure~\ref{fig:1} on the right. Let $\sigma_{|{\mathrm{bottom}}}$ and $\sigma_{|{\mathrm{top}}}$ be the restrictions of $\sigma$ to the top and bottom of $\mathbb S(N,M)$, considered as two elements of $\{-1,1\}^{(\mathbb Z/N\mathbb Z)^{d-1}}$.
Introduce the quantity
\begin{equation}
Z(N,M,\underline\tau,\overline\tau):=\sum_{\sigma\in \{-1,1\}^{\mathbb S(N,M)}}\exp[-\beta H_{\mathbb S(N,M),h}(\sigma)]\mathbb I(\sigma_{|\mathrm{bottom}}=\underline\tau,\sigma_{|\mathrm{top}}=\overline\tau),
\end{equation}
where $\underline\tau,\overline\tau\in\{-1,1\}^{(\mathbb Z/N\mathbb Z)^{d-1}}$, as well as the so-called {\em transfer matrix }
\begin{equation}
V_N(\underline\tau,\overline\tau):=Z(N,1,\underline\tau,\overline\tau)=\exp\Big[-\beta\Big( \sum_{x\in (\mathbb Z/N\mathbb Z)^{d-1}}\underline\tau_x\overline\tau_x+ \sum_{\{x,y\}\in E((\mathbb Z/N\mathbb Z)^{d-1})}\overline\tau_x\overline\tau_y\Big)\Big].
\end{equation}
One immediately finds that 
$
Z(N,M,\underline\tau,\overline\tau)=V_N^M(\underline\tau,\overline\tau).
$
Other quantities of the model may be written in terms of transfer matrices: for instance the partition function of the model on the $d$-dimensional torus $(\mathbb Z/N\mathbb Z)^d$ becomes the trace of $V_N^N$. 

One important aspect of those transfer matrices $V_N$ is that certain questions on the behaviour of the model are rephrased as spectral questions on the transfer matrix. For instance, by letting $M$ and then $N$ go to infinity, one observes that the asymptotic behaviour of the partition function on $(\mathbb Z/N\mathbb Z)^{d-1}\times(\mathbb Z/M\mathbb Z)$, and therefore the value of the free energy, are connected to the asymptotic behaviour of the leading eigenvalue of $V_N$ as $N$ tends to infinity. This can very well be an intractable problem, but in some cases it is not.  
  \end{framed}
  
\subsection{Onsager's result}\label{sec:3.2}

Kramers and Wannier's results unraveled the potential mathematical interest of the Ising model, but the real revolution came only a few years after with one of the most impressive achievements in mathematical physics.
Lars Onsager, Nobel prize winner in 1968, was a Norwegian specialist in theoretical chemistry. He was particularly interested in mathematical problems and focused his attention on the Ising model for the formidable challenge that its exact solution represented more than for his physical relevance\footnote{His opinion on this fact seemingly evolved: in his first 1944 paper \cite{Ons44} he presents it as a poor model of ferromagnetism, but a fairly good model for binary alloys, while in his paper with Kaufman in 1949 \cite{KauOns49} he describes it as a model for  ferromagnetism.}. 

To everyone's surprise, Onsager announced at a conference of the New York Academy of Sciences in 1942 that he obtained the following {\em exact} expression for the free energy (at zero magnetic field) of the Ising model on the square lattice $\mathbb Z^2$:
\begin{equation}
-\beta f(\beta,0)=\ln 2+\frac1{8\pi^2}\int_0^{2\pi}\int_0^{2\pi}\ln[\cosh(2\beta)^2-\sinh(2\beta)(\cos\theta_1+\cos\theta_2)]d\theta_1d\theta_2.
\end{equation}
This implies, in physics jargon, that the model is {\em exactly solvable}. This solvability is itself linked to a deep property of the model called {\em integrability}.
Onsager's underlying idea was that the transfer matrices of the 2D Ising model are a product of two matrices that generate (by taking successive brackets) a finite-dimensional Lie algebra. He used this observation to derive the asymptotic behaviour of the leading eigenvalue of these matrices in his famous 1944 paper \cite{Ons44}. In 1949, Bruria Kaufman \cite{KauOns49} provided an alternative and simpler derivation.

A few years later, Onsager surprised the world of theoretical physicists again by claiming an exact expression for the spontaneous magnetisation on $\mathbb Z^2$: for $\beta\ge \beta_c$,
\begin{equation}
m^*(\beta)=(1-\sinh(2\beta)^{-4})^{1/8}.
\end{equation}
While the result was announced by Onsager first, it was a young physicist, that would later become one of the most influential theoretical physicists of the second half of the twentieth century, Chen-Ning Yang (from Yang-Baxter's equation, Yang-Mills's theory, Lee-Yang's theory, etc.), who provided a mathematical proof \cite{Yan52} of this statement by achieving a mathematical tour de force involving Toeplitz determinants. The proof relies on a computation, again using transfer matrices but much more evolved than for the free energy, of the two-point function $\langle\sigma_{(0,0)}\sigma_{(n,0)}\rangle_{(\mathbb Z/N\mathbb Z)^2,\beta,0}$, and the observation that 
\begin{equation}
m^*(\beta)^2=\lim_{n\rightarrow\infty}\lim_{N\rightarrow\infty}\langle\sigma_{(0,0)}\sigma_{(n,0)}\rangle_{(\mathbb Z/N\mathbb Z)^2,\beta,0}
\end{equation}
(at the time, such an identity was not obviously true, but nowadays this can be proved easily using for instance the FK percolation, see Section~\ref{sec:7.2}).

In the forties and fifties, these successes were considered by physicists as a mathematical curiosity rather than a truly crucial advance. Yet, they had a revolutionary impact on theoretical physics for multiple reasons:
First, the {\em level of sophistication} of the mathematical tools used in the proofs is without any common measure with what was previously used in such kinds of problems, and these techniques created whole new types of mathematical physics. 
Second, the behaviour of the model did not correspond to previous {\em mean-field approximations}, thus invalidating rigorously the Curie-Weiss  or Landau theories and opening a new era in statistical mechanics. 
 Third, the results had {\em many direct applications} for the Ising model itself: for instance the specific heat $C(\beta)$ can easily be shown to blow up logarithmically as $\beta$ approaches $\tfrac12\log(1+\sqrt 2)$, thus confirming rigorously that this value is the critical point of the system (the logarithmic blow up is one example of non mean-field behaviour).

Numerous alternatives have been proposed to the approach of Onsager-Kaufman-Yang, often referred to as the {\em algebraic method}. As a joke, Baxter and Enting named their 1978 paper \cite{BaxEnt78}, introducing a solution to the 2D Ising model involving the notion of star-triangle transformation, the "{\em 399th solution of the Ising model}". This count is, of course, overestimated, but one can list a large number of alternative strategies.

The first such strategy is called the {\em combinatorial approach} and is referring to an original argument of Kac and Ward \cite{KacWar52} rewriting the partition function of the model in terms of the square root of the determinants of so-called Kac-Ward's matrices using a combinatorial expansion of the partition function generalising the van der Waerden high-temperature expansion \cite{Wae41}. The advantage of such an approach is that it does not rely on transfer matrices, and therefore is applicable to {\em every} finite planar graph, even with arbitrary nearest-neighbour coupling constants. Unfortunately, the original argument was not entirely rigorous and one had to wait until 1999 \cite{Dob99} to finally obtain a mathematical derivation of this approach. Nowadays, the method is very well understood and especially useful in relation to discrete holomorphicity and higher genus graphs, see \cite{Cim12} and references therein for a more complete account.

The (nowadays) most classical method is probably the {\em Pfaffian method}. It came as an attempt to go around the substantial difficulties to make the combinatorial approach rigorous. Due to Hurst and Green \cite{HurGre60}, Kasteleyn \cite{Kas63}, and Fisher \cite{Fis67}, the strategy consists in writing the Ising partition function on a finite planar graph $G$ in terms of the dimer (a dimer configuration is a subset of edges which covers every vertex exactly once) partition function on a related graph $K(G)$ (the precise definition of the graph depends on the implementation of the Pfaffian method). It is then possible to relate the partition function to a skew-symmetric adjacency matrix and express the partition as a Pfaffian, hence the name of the method. This strategy has been the basis of a number of more refined results about the model, in particular the computation of the spin-spin correlations of the model at and away from criticality. For the deepest and most impressive results, we recommend that the reader takes a look at the two books of McCoy-Wu \cite{McCWu13} and Palmer \cite{Pal07}.

Another approach of importance was proposed by Schultz-Mattis-Lieb in \cite{SchMatLie64} to tackle the cases for which a transfer matrix can be used. In this paper, they connected the transfer matrix with the exponential of a quantum hamiltonian. This connection to 1D quantum spin chains has been very fruitful and understood in a number of alternative ways since then. As a byproduct, the authors were able to express the partition function as a Grassmann "Gaussian" integral. The advantage of this way of writing the partition function is that the Pfaffians emerge naturally. This approach is at the basis of renormalisation schemes in two dimensions; see Section~\ref{sec:8.3}.

Yet another approach dealing with the context in which transfer matrices can be applied is worth mentioning, as it is by far the most generalisable to other models. It is based on the commutation of the transfer matrices attached to the model with different critical parameters. Pioneered by Rodney Baxter, this approach consists in using the so-called Yang-Baxter equation. The advantage is that the same strategy can be applied to a very large variety of {\em integrable systems}, such as the six-vertex model, to cite only one example. We refer to \cite{Bax89} and references therein for more details.

\section{The fifties and sixties: the Ising model becomes a laboratory for understanding critical phenomena}\label{sec:4}
The fifties and sixties were probably the decades during which the Ising model became an "unavoidable" model. The realisation that having a tractable model of statistical physics could be a useful explanatory but also predicting tool became more and more obvious. The Ising model, with Onsager's solution, was a prime example of a model with such qualities.

The model therefore developed tremendously in the postwar era in theoretical physics as well as in a rapidly growing field called {\em mathematical physics}. The latter gathered more and more physicists that were interested in rigorous aspects of the objects they studied, and mathematicians willing to study problems that were naturally emerging from physical modelling. The Ising model offered a wonderful playground for such scientists, and the number of papers mentioning the model started to be counted in the hundreds. 

\subsection{Progress in mathematical physics:  from perturbative regions of the phase diagram to the vicinity of the critical point}\label{sec:4.1}

During this period, the newly developing community of mathematical physicists recognised that the study of phase transitions, and in particular of the critical phase (when $\beta$ is equal to $\beta_c$), was a vast field of its own. While the previous developments mostly concerned values of $\beta$ and $h$ that were far from the critical regime (Peierls' argument \cite{Pei36} or Baker's use of Pad\'e approximant \cite{Bak61} for instance), the situation changed drastically around the fifties. The interest in the intermediate values of $\beta$ became stronger and stronger. 
Onsager's solution offers a precise understanding of the critical behaviour of the 2D Ising model, yet it has clear downsides related to the relative fragility of the integrability of the system. As a consequence, mathematical physicists started using the Ising model not only as a solvable system, but more generally as a good mathematical model that one should not reduce to its integrability aspects. New rigorous techniques emerged during this period to try to understand the vicinity of the critical point for non-integrable cases, for instance in higher dimensions.

\subsubsection{Correlation inequalities}\label{sec:4.1.1}

It is natural to ask which monotonicity properties are satisfied by the system, in particular by the {\em spin-spin correlations} 
$\langle\sigma_A\rangle_{G,\beta,h}$ where $\sigma_A:=\prod_{x\in A}\sigma_x,
$
when the parameters vary (for instance, $G$, $\beta$ or $h$). 

To tackle such questions, mathematical physicists started proving  what we now call {\em correlation inequalities} using combinatorial arguments. Among the first such examples are Griffiths' inequalities \cite{Gri67}: for every $\beta,h\ge0$ and every $A,B\subset V$,
\begin{equation}
\langle\sigma_A\rangle_{G,\beta,h}\ge0\qquad\text{ and }\qquad\langle\sigma_A\sigma_B\rangle_{G,\beta,h}\ge \langle\sigma_A\rangle_{G,\beta,h}\langle\sigma_B\rangle_{G,\beta,h}.
\end{equation}
A byproduct of the second inequality when applied to $B=\{x,y\}$ and summed over all edges $\{x,y\}$, is that correlations $\langle\sigma_A\rangle_{G,\beta,h}$ are increasing in $\beta$ (and also in $G$ with a little bit of additional work). One may derive the same for the spontaneous magnetisation $m^*(\beta)$, so that the definition of $\beta_c$ can now be rephrased as
\begin{equation}
\beta_c=\inf\{\beta\ge0:m^*(\beta)>0\}=\sup\{\beta\ge0:m^*(\beta)=0\}.
\end{equation}
This implies in particular that there is indeed a unique transition between paramagnetic and ferromagnetic phases.

Other interesting correlation inequalities were obtained in subsequent years. Let us contemplate a few examples (we do not write them in full generality, and we drop the subscript after $\langle\cdot\rangle$):
\begin{itemize}
\item GHS's inequality \cite{GriHurShe70}: for $h\ge0$ and $x\in G$,
\begin{equation}
\tfrac{\partial^2}{(\partial h)^2}\langle\sigma_x\rangle\le 0.
\end{equation}
\item Simon-Lieb's inequality \cite{Lie80}: for $S\ni0$ and $x\notin S$, when $\langle\cdot\rangle_S$ refers to the model in $S$,
\begin{equation}
\langle\sigma_0\sigma_x\rangle\le \sum_{y\in \partial S}\langle\sigma_0\sigma_y\rangle_S\langle\sigma_y\sigma_x\rangle.
\end{equation}
\item Messager-Miracle-Sol\'e's inequality \cite{MesMir77}: for  $x,y\in \mathbb Z_+^d$ (below $\langle\cdot\rangle$ is defined on $\mathbb Z^d$)
\begin{equation}
\langle\sigma_0\sigma_{x+y}\rangle\le \langle\sigma_0\sigma_x\rangle.
\end{equation} 
\item FKG's inequality \cite{ForKas72}: for  any increasing functions $f,g:\{-1,1\}^V\rightarrow \mathbb R$,
\begin{equation}
\langle fg\rangle\ge \langle f\rangle\langle g\rangle.
\end{equation} 
\end{itemize}
This far from exhaustive list, which we did not discuss in detail, is intended to show the variety of possible correlation inequalities.
Clever use of these inequalities provided the embryo of what would be considered later as the theory of non-critical statistical physics systems at equilibrium, as the correlation inequalities and their consequences often generalise in the same (or slightly altered) form  to a wider class of lattice spin models.
\subsubsection{The Ising model with a magnetic field: the Lee-Yang theory}\label{sec:4.1.2}

While studying the whole phase diagram is a Herculean task that was far beyond the techniques developed at the time, a beautiful development enabled mathematical physicists to understand the case $h\ne 0$. 

The twin papers \cite{LeeYan52}, referred to as the Lee-Yang theory,  relate the regularity properties of the free energy (and therefore the location of singular points corresponding to places where a phase transition occurs) to the locus of the complex zeroes of the partition function $Z(G,\beta,h)$ when seen as a function of $h\in\mathbb C$. Beyond the result itself, the philosophy consisting in studying the complex zeroes of the partition function had a resounding effect on the field of mathematical physics. This can be put in parallel with the analysis of zeroes of the Riemann zeta function: one learns something about prime numbers by studying the zeroes of a generating-type function associated with them. 

The result of Lee and Yang is not restricted to the n.n.f.~Ising model on $G\subset\mathbb Z^d$, but the latter gives an important application of it. In our context,  let $Z(G,\beta,\mathbf h)\in\mathbb C$ (for $\mathbf h=(\mathbf h_x:x\in V)\in\mathbb C^V$) be the partition function defined as in \eqref{eq:def2} with the difference that the magnetic field is allowed to vary with the vertex, i.e.~that the $\sum_{x\in V}h\sigma_x$ term of the Hamiltonian in \eqref{eq:def1} is replaced by $\sum_{x\in V}\mathbf h_x\sigma_x$. 
The result 
states that for this model, the zeroes of the function
$\mathbf h\mapsto Z(G,\beta,\mathbf h)$ are satisfying $\mathrm{Re}(\mathbf h_x)= 0$ for every $x\in V$.

As a consequence of this theorem, the free energy $f(\beta,h)$ (which we recall from \eqref{eq:free energy} is expressed in terms of the limit of the logarithm of partition functions) is analytic as soon as $h\ne0$. Other consequences follow, such as exponential decay of so-called truncated correlations of the system, as well as analyticity of the other thermodynamical quantities when the magnetic field is non-zero. Roughly put, the Lee-Yang theory enables to understand in full detail  the part of the phase diagram corresponding to non-zero magnetic field.

\subsection{Revolutionary progress on the physics front}\label{sec:4.2}

In parallel to these first successes in mathematical physics, revolutionary progress was made during this period on the physical understanding of phase transitions. Among other things, the scaling and universality hypotheses were formulated, and the pillars of the renormalisation group were cast, in both cases using the Ising model as an important source of inspiration.

\subsubsection{Critical exponents and the success of scaling theory}\label{sec:4.2.1}

A fundamental notion of physics is the assumption that thermodynamical quantities of physical systems near criticality tend to take simple forms when expressed in terms of the parameters of the system. A major advance was achieved in the sixties by American chemist Benjamin Widom who proposed in \cite{Wid65} that these quantities are powers in each parameter. For the Ising model, the parameters are $\beta$ and $h$, and this {\em scaling hypothesis} translates into the existence of so-called {\em critical exponents}. To give a few examples related to already defined quantities, one may for instance predict that \begin{equation}\label{eq:expo}
m^*(\beta)=(\beta-\beta_c)_+^{\bm\beta+o(1)}\quad\quad m(\beta_c,h)=h^{1/\bm\delta+o(1)}\quad\quad\langle\sigma_0\sigma_x\rangle_{\beta_c,0}=\frac{1}{|x|^{d-2+\bm\delta}}
\end{equation}
(notice that $\bm\beta$ and $\beta$ have nothing to do with each other), where $o(1)$ is a quantity tending to 0 as $\beta$ tends to $\beta_c$, $h$ tends to 0, or $|x|$ tends to infinity, respectively. In fact, a whole family of such exponents, denoted by $\bm\alpha,\bm\beta,\bm\gamma,\bm\delta,\bm\eta,\bm\nu$ (for the most classical ones) can be defined for each model.
Understanding the phase transition boils down to, among other things,  deriving those exponents.

Dealing with such exponents, one may naturally wonder how many degrees of freedom truly exist in statistical physics models. For instance, could some of these critical exponents be connected via direct relations that would transcend the precise definition of each model? In the sixties, physicists such as Essam, Fisher, and Widom himself, to cite only those three (see \cite{EssFis63,Fis64,Wid65} for some early works on the subject), started unraveling systematic connections between the exponents, thus hinting towards the fact that only two degrees of freedom exist and that exponents are related by so-called {\em scaling relations}
\begin{equation}
\bm\nu d=2-\bm\alpha=2\bm\beta+\bm\gamma=\bm\beta(\bm\delta+1)=\bm\gamma\frac{\bm\delta+1}{\bm\delta-1}\qquad2-\bm\eta=\frac{\bm\gamma}{\bm\nu}=d\frac{\bm\delta-1}{\bm\delta+1}.
\end{equation}
The scaling relations apply in a context which is far more general than just the Ising model (see for instance \cite{DumMan20} for a proof in the case of a large family of two-dimensional percolation models). In the course of discovering these different scaling relations, the Ising model in two and three dimensions played the important role of a sanity check. While other experimental systems were used as testing grounds, the Ising model was the only example of a theoretical system which did not exhibit mean-field behaviour (and therefore was not too ``trivial'') and for which such exponents were available, either rigorously thanks to the exact solution in 2D, or approximately thanks to Baker's use of Pad\'e approximant \cite{Bak61} in 3D.

To conclude this section, let us mention an important quantity, called the {\em correlation length} $\xi(\beta)$ of the system, that plays an important role in the scaling hypothesis (it corresponds to the exponent $\bm\nu$). 
We consider the case $\beta<\beta_c$ but a similar notion can be introduced for $\beta>\beta_c$, with analogous interpretations.

When considering, say, spin-spin correlations at criticality, one expects an algebraic decay as mentioned in \eqref{eq:expo}. Yet, when $\beta<\beta_c$, the scaling hypothesis cannot hold uniformly in $|x|$ and such a decay does not occur. In fact, it was found in many systems that spin-spin correlations decay exponentially fast  (see Section~\ref{sec:7.1} for more details) and the inverse-rate of decay is the correlation length $\xi(\beta)$. 
This correlation length has an interesting interpretation: it is the smallest scale at which the system with $\beta<\beta_c$ is off-critical, meaning that when looking at a system with a size which is much smaller than $\xi(\beta)$, the difference between the system and a critical system will be invisible to the physicist's eye, while on the contrary when the size is much larger than $\xi(\beta)$, the model looks similar to the case of $\beta\ll \beta_c$. In other words, when approaching the critical point, a system becomes more and more ``critical''. By how much this is true depends on the size of the system, and the correlation length separates between the sizes at which the system looks critical, and the sizes at which it looks clearly non critical.

\subsubsection{Kadanoff's block-spin renormalisation and universality}\label{sec:4.2.2}

While Widom's scaling hypothesis provides compelling evidence that critical exponents exist, the underlying justification of the hypothesis itself remained slightly superficial until Russian physicist Leo Kadanoff provided an illuminating argument for it. In his famous 1966 paper \cite{Kad66}, Kadanoff suggested that the block-spin renormalisation transformation -- i.e.~replacing a block of neighbouring sites by one site having a spin equal to the dominant spin in the block -- corresponds to appropriately changing the scale and the 
parameters $\beta$ and $h$ of the model. Assuming that iterating this procedure somehow converges suggests that the asymptotic properties of the system are described by a fixed point of a renormalisation map. As a result, one ends up with the {\em scale invariance} of the model. This argument, inspired by the study of the Ising model, turned out to be the basis of the monumental theory of the renormalisation group (RG) that was put in a general framework a few years later by Kenneth Wilson \cite{Wil71}. 

The block-spin argument of Kadanoff achieved much more than a physical justification of the scaling hypothesis. Assuming {\em uniqueness of the fixed point} also implies that the renormalisation of Ising models defined on different $d$-dimensional lattices should converge to the same fixed point, and therefore share the same critical exponents. This was already partially realised in 2D by observing the Ising model on the square, hexagonal, and triangular lattices (they are all exactly solvable) as well as in 3D by approximations using series expansions \cite{DomSyk57}, but the renormalisation argument suggests that the few examples of equalities between exponents are, in fact, the illustration of a much more general phenomenon. 

What is now known as the {\em universality hypothesis} was explicitly formulated in parallel by Robert B. Griffiths and Kadanoff in 1971 \cite{Gri70,Kad71}. Roughly speaking, it states that the critical properties of a physical system only depend on 
\begin{itemize}[noitemsep]
\item the lattice dimension $d$;
\item the symmetry of the space of possible spins ($\mathbb Z/2\mathbb Z$ symmetry for Ising);
\item the speed of decay of coupling constants (this is only relevant when the $J_{x,y}$ are allowed to decay polynomially with $\|x-y\|$, which is not the case in this text).
\end{itemize}
This realisation of universality is fundamental to the relevance of statistical physics as a whole. To borrow from Kadanoff's wording: ``{\em Why study a simplified model like the Ising model? The strategy of studying physical questions by using highly simplified models is made rewarding by a characteristic of physical systems called ``universality'', in that many systems may show the very same qualitative features, and sometimes even the same quantitative ones. To study a given qualitative feature, it often pays to look for the simplest possible example.}''

To summarise Section~\ref{sec:4}, by the end of the sixties it became clear to mathematical physicists and theoretical physicists that the Ising model was one of the most striking examples of a simple physical system which was rich enough to grasp a large variety of phenomena falling in the range of statistical physics. Results on the Ising model started to play a role similar to experimental results in the sense that they could corroborate or, on the contrary, invalidate the embryo of a theory.
It is fair to say that the importance of the model was never argued upon later on and that it was finally recognised as one of the centerpieces of modern statistical physics. 

\section{The sixties and seventies: emergence of the probabilistic interpretation}\label{sec:5}

Physicists and mathematical physicists think of the quantity $\langle\cdot\rangle_{G,\beta,h}$ as a form attributing to each function $X:\{-1,1\}^V\rightarrow\mathbb R$ (resp.~$\mathbb C$) a value in $\mathbb R$ (resp.~$\mathbb C$). In the late sixties and seventies, the rise of probabilistic methods led to an alternative interpretation of the Ising model in which $\langle\cdot\rangle_{G,\beta,h}$ is now understood as (dual to) a probability measure $\mu_{G,\beta,h}$. As a consequence of this reinterpretation, it becomes natural to ask what the properties of a randomly chosen spin configuration are, and what the possible measures on the infinite lattice that can be obtained as limits of measures in finite volume are.

\subsection{The random geometry of the spin configuration} \label{sec:5.1}

As mentioned above, $\langle\cdot\rangle_{G,\beta,h}$ is the linear form associated with the probability measure $\mu_{G,\beta,h}$ on $\{-1,1\}^V$ defined for every configuration $\sigma$ by the formula
\begin{equation}
\mu_{G,\beta,h}[\{\sigma\}]:=\frac{1}{Z(G,\beta,h)}\exp[-\beta H_{G,h}(\sigma)].
\end{equation}
Then, quantities like $\langle\sigma_A\rangle_{G,\beta,h}$ can be interpreted as the {\em correlations} between the random variables $\sigma_x$ with $x\in A$. Note that in this interpretation the partition function is a normalising factor making the measure at hand a probability measure.

Let us assume for a moment that $h=0$ and interpret the phase transition in terms of probability. The structure of the probability measure is such that configurations have greater probability if they have more pairs of neighbours with a similar spin. In this interpretation, the larger $\beta$ is the more important it is that neighbours have the same spins. In particular, in the limit as $\beta$ tends to infinity, one ends up with one of the two configurations where all spins are the same. It becomes then natural to expect that for $\beta$ large, typical configurations have an excess of one spin compared to the other. On the other hand, when $\beta$ is very small, how much the measure takes the agreements into account is fairly limited, and one may expect that spins behave roughly independently, at least at large distance of each other. 

The interpretation in terms of random variables opens new uncharted territories: one can interpret probabilistically natural thermodynamical quantities such as magnetisation (which corresponds to the expectation of the spin at a vertex) or surface tension. It also opens a way to new problems, such as dynamics on the space of spin configurations or  large deviations (for instance for an Ising model at an inverse-temperature $\beta$, but with an excess of $+1$ spins in a region and of $-1$ spins in another); see Frame 3.

\begin{figure}[t]
  \begin{center} \includegraphics[width=0.65
  \textwidth]{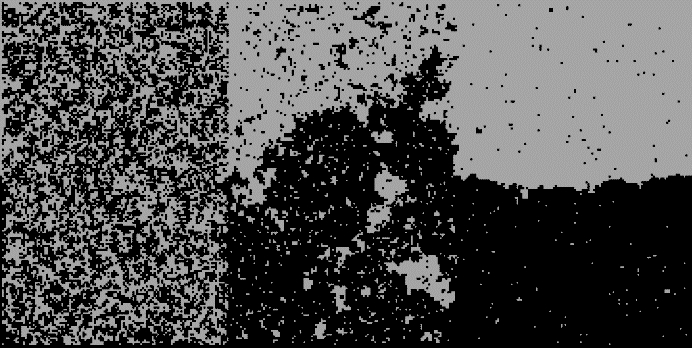}\ \includegraphics[width=.327\textwidth]{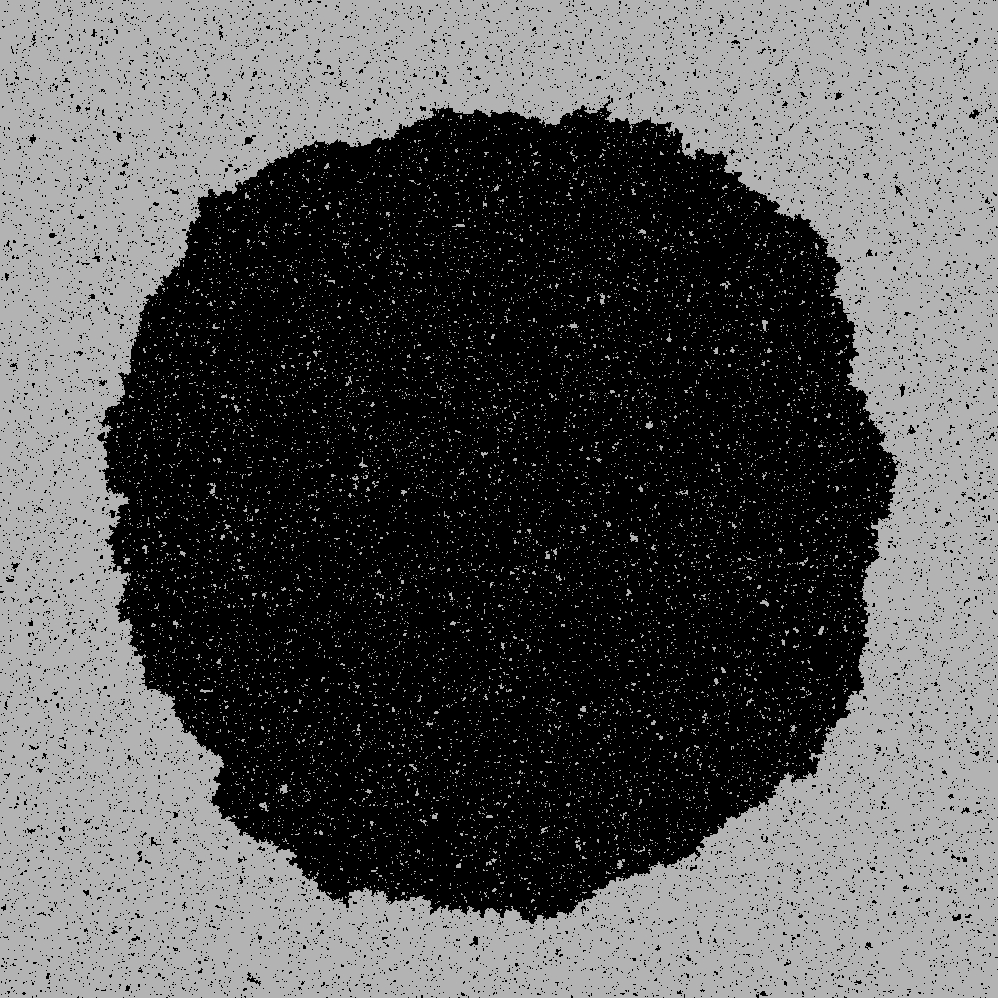}\caption{\textbf{On the left.} Simulations at three different temperatures ($\beta<\beta_c$, $\beta=\beta_c$, and $\beta>\beta_c$) of the Ising model with plus boundary conditions on the top and minus boundary conditions on the bottom. Pluses are in gray and minuses in black. {\em Credit: S.~Smirnov}. \textbf{On the right.} An example of a bubble of minuses in an environment of pluses at $\beta>\beta_c$. {\em Credit: Y.~Velenik.}}\label{fig:2}\end{center}
 \end{figure}

    \begin{framed}
  \centerline{\textbf{Frame 3: sampling the Ising model -- Glauber dynamics}
}
The probabilistic interpretation naturally raises the question of sampling random configurations according to $\mu_{G,\beta,0}$ (set $h=0$ for simplicity). A classical method consists in expressing the measure as the invariant measure of a Markovian dynamics $(\sigma(t):t\ge0)\in (\{-1,1\}^V)^{\mathbb R_+}$, called the Glauber dynamics and defined as follows: attach an exponential clock to each vertex of $G$. Each time a clock rings, say at time $t$ at $x\in V$,
\begin{itemize}[noitemsep]
\item If $\sigma_x(t)\sum_{y:\{x,y\}\in E}\sigma_y(t)<0$, switch the value of the spin at $x$,
\item Otherwise, switch the value of the spin at $x$ with a probability equal to $\exp[-2\beta \sum_{y:\{x,y\}\in E}\sigma_y(t)]$, and do not switch otherwise.
\end{itemize}
Since $\mu_{G,\beta,0}$ is the only invariant measure for this dynamics, the limit as $t$ tends to infinity, irrespectively of the initial  value $\sigma(0)$, is sampled according to $\mu_{G,\beta,h}$.

This dynamics was named after the American physicist Roy J.~Glauber. Alternative choices of dynamics are obtained by changing the jump probabilities. In Figure~\ref{fig:2}, three simulations of the Ising model are shown respectively below (on the left), at (in the middle) and above (on the right) $\beta_c$.
  \end{framed}

 \subsection{Boundary conditions and the Gibbs formalism}\label{sec:5.2}

An important output of the probabilistic interpretation of the model is that it becomes natural to {\em condition} on spins in a subset of $V$. More precisely, let $W\subset V$ and let $H$ be the graph with vertex-set $W$ and edge-set induced by the edges of the graph $G$. Let $\tau\in\{-1,1\}^V$ be a spin configuration on $G$. One may ask what is the law of the spins in $W$ when conditioning  $\sigma$ outside $W$ to be equal to $\tau$, i.e.~what is $\mu_{G,\beta,h}[\,\cdot\,|\sigma_x=\tau_x,\forall x\notin W]$?

The answer to this question is best cast when introducing the notion of boundary conditions. For a subgraph $G$ of $\mathbb Z^d$ and a configuration $\tau\in\{-1,1\}^{\mathbb Z^d}$, introduce the measure $\mu_{G,\beta,h}^\tau$ with $\tau$ {\em boundary conditions} defined like $\mu_{G,\beta,h}$ except that $H_{G,h}$ is replaced by 
\begin{equation}
H_{G,h}^\tau(\sigma):=H_{G,h}(\sigma)-\sum_{\{x,y\}\in E(\mathbb Z^d): x\in V, y\notin V}\sigma_x\tau_y.
\end{equation}
Note that the only values of $\tau$ that matter are on the exterior boundary of $G$, i.e.~on the vertices that are connected by an edge of $\mathbb Z^d$ to a vertex in $V$.

With this definition, we obtain the following important property of the Ising model, called the {\em spatial Markov property}: 
for every finite subgraph $G$ of $\mathbb Z^d$, every $W\subset V$, and every configuration $\tau\in\{-1,1\}^{\mathbb Z^d}$, if $H$ denotes the graph induced by the set $W$,
\begin{equation}
\mu_{G,\beta,h}[\,\cdot\,|\sigma_x=\tau_x,\forall x\notin W]=\mu_{H,\beta,h}^\tau[\,\cdot\,].
\end{equation}
In words, when conditioning the Ising model on $G$ to coincide with a given configuration
outside $W$, one gets the measure in $H$ with the corresponding boundary
condition.

This property offers a natural consistency relation between measures $\mu_{G,\beta,h}^\tau$ for varying $\tau$ and $G$. As a byproduct, one is naturally led to postulate that any reasonable infinite-volume version of Ising measures should satisfy the same consistency relation. One therefore ends up with the following notion: a measure $\mu$ on $(\{-1,1\}^{\mathbb Z^d},\mathfrak F_{\mathbb Z^d})$ is called a {\em Gibbs measure} of the Ising model with parameters $\beta$ and $h$ if it satisfies the {\em Dobrushin-Lanford-Ruelle (DLR) property}: for every finite $V\subset\mathbb Z^d$ and $\tau\in\{-1,1\}^{\mathbb Z^d}$,
\begin{equation}
\mu[\,\cdot\,|\mathfrak F_{\mathbb Z^d\setminus V}]=\mu_{G,\beta,h}^{\tau}[\,\cdot\,]\text{ on }E_\tau\text{ $\mu$-almost surely},
\end{equation}
where 
\begin{itemize}[noitemsep,nolistsep]
\item $G$ is the graph induced by  the vertex-set $V$; 
\item $E_\tau$ is the event that $\sigma$ and $\tau$ agree on the exterior boundary of $G$; 
\item $\mathfrak F_{\mathbb Z^d\setminus V}$ is the $\sigma$-algebra generated by the random variables $(\sigma_x:x\notin V)$.\end{itemize}
The notion of Gibbs measure is not restricted to the Ising model (see \cite{Geo11} for a book on the subject), but the classification of such Gibbs measures has been the object of intense study in the specific case of the Ising model, with a very successful outcome. 

The first question that one may ask is the existence of Gibbs measures. At least three such measures can be defined in a fairly straightforward way. By taking limits as $G$ tends to $\mathbb Z^d$ of the measures $\mu_{G,\beta,h}$, $\mu_{G,\beta,h}^+$, and $\mu_{G,\beta,h}^-$ (where $+$ and $-$ refer, with a slight abuse of notation, to $\tau$ equal to all $+1$ or all $-1$), one ends up with three (possibly equal) Gibbs measures $\mu_{\beta,h}$, $\mu_{\beta,h}^+$ and $\mu_{\beta,h}^-$. More generally, one may construct measures by taking all possible sub-sequential limits of measures of the form $\mu_{G,\beta,h}^\tau$, where one may even consider $\tau$ as a random variable. 

In general, the set of possible Gibbs measures on $\mathbb Z^d$ is a non-empty simplex whose extremal measures are called {\em extremal states}. One can therefore try to classify such extremal Gibbs measures.

Some cases are quite simple to treat: for $h\ne0$ or $h=0$ and $\beta<\beta_c$ the simplex is reduced to a singleton, i.e.~there exists a {\em unique} Gibbs measure. When $h=0$ and $\beta=\beta_c$, it was recently proved that this is also the case \cite{AizDumSid15}. On the contrary, when $h=0$ and $\beta>\beta_c$, things are more interesting. It was realised very early on that there may be more extremal states than the two obvious $\mu_{\beta,0}^+$ and $\mu_{\beta,0}^-$, but examples that were found did not exhibit translation invariance. The most important such specimen was provided by Russian mathematical physicist Roland Dobrushin \cite{Dob73}, who explained that in three dimensions, the measure $\mu_{\beta,0}^{\mathrm{dobr}}$ obtained by taking the limit of measures $\mu_{[-n,n]^3,\beta,0}^\tau$, where $\tau$ is all plus on the upper half-space, and all minus on the lower half-space, was not translation invariant in the vertical direction at high values of $\beta$. The existence of these Dobrushin states is related to a very deep and still mysterious (at least on a mathematical level) phenomenon in 3D statistical physics often referred to as the {\em roughening phase transition}.

Leaving non-translation invariant measures aside, many efforts were made to prove that every {\em translation invariant} Gibbs state is a convex combination of $\mu_{\beta,0}^+$ and $\mu_{\beta,0}^-$. The first result in this direction proved a stronger statement that draws a direct link with the previous paragraph. In two dimensions, Aizenman \cite{Aiz80} and Higuchi \cite{Hig78} proved in the eighties that {\em every} Gibbs state, not only translation invariant ones, is a mixture of $\mu_{\beta,0}^+$ and $\mu_{\beta,0}^-$. In particular, $\mu_{\beta,0}=\tfrac12\mu_{\beta,0}^++\tfrac12\mu_{\beta,0}^-$. In higher dimensions, it took twenty more years to obtain the result for every {\em translation invariant} Gibbs measure. We refer to the historical proof of Bodineau \cite{Bod06} and to the recent generalisation of Raoufi \cite{Rao20}.

\subsection{Phase coexistence and Wulff shape}\label{sec:5.3}

The classification of Gibbs states naturally raises the question of the coexistence of different so-called {\em phases}. When $h=0$ and $\beta>\beta_c$, $\mu_{\beta,0}^+$ and $\mu_{\beta,0}^-$ are not equal: they correspond to two extremal states, sometimes referred to as the plus and minus phases. Now, what happens when one tries to ``mix'' the two states? For instance, how does it look if one asks that part of the space is in one state, and the other part is in the other one? 

In 2D, an {\em interface} is created between the two phases (see Figure~\ref{fig:2} for simulations at different temperatures). While it is not obvious to define such an object in general, let us consider the simple example of the Ising model on a finite box $[-n,n]^2$ of the triangular lattice with plus spins on the part of the boundary above the $x$-axis, and minus spins on the rest of the boundary. In this case, one can draw a unique interface going from $(-n,0)$ to $(n,0)$ winding between pluses and minuses. It was understood  heuristically early on that  above criticality  this interface should have the same fluctuations as Brownian motions, but it took decades to turn this intuition into a rigorous proof, first in the large $\beta$ regime and then in the whole $\beta>\beta_c$ regime; see \cite{GreIof05} and references therein. The techniques involved also enabled mathematicians to understand precise asymptotics of spin-spin correlations in the non-critical regimes. The theory, known under the coined name of Ornstein-Zernike theory, is now an area of intense research and spans over a large variety of statistical physics models. We refer to \cite{CamIofVel03} for details on the Ising case.

When conditioning on the neighbourhood of the origin to be in a plus phase inside a minus phase, one ends up with a ``bubble'' (see Figure~\ref{fig:2} on the right) converging when taking larger and larger volume to the so-called {\em Wulff shape}. In 2D, this bubble was analysed in detail, see the book \cite{DobKotShl92} and the article \cite{IofSch98}.
In 3D, the story is even more complex. The boundary between the plus and minus phases is a kind of two-dimensional surface.  The study of this object is quite intricate, and the fluctuations of the surface are still widely open. We refer to \cite{Bod99,CerPiz00,BodIofVel00} and references therein.

\section{The seventies and eighties: the Ising model and field theory}\label{sec:6}

\subsection{Constructive quantum field theory}\label{sec:6.1}

Quantum field theories with local interaction are central in most  subfields of theoretical physics, from high energy to condensed matter physics.  
The mathematical challenge of the proper formulation of this concept led to the program of constructive quantum field theory (CQFT).  A path towards that goal was charted through the  proposal to define quantum fields satisfying Wightman axioms~\cite{Wig56} using the Osterwalder-Schrader theorem~\cite{OS73}, in which case the construction boils down to producing relevant random distributions defined over the corresponding Euclidean space that meet a number of conditions such as suitable analyticity, permutation symmetry, Euclidean covariance, and reflection-positivity.    

Finding these \emph{Euclidean fields} boils down to constructing probability averages over random distributions $\Phi(x)$ of the form\begin{equation} \label{Phi_EV}
\langle F(\Phi) \rangle  \approx \frac1{\mathrm{norm}}   \int  F(\Phi) \exp[-H(\Phi)]  \prod_{x\in \mathbb R^d} d \Phi(x) ,
\end{equation} 
where 
\begin{itemize}[noitemsep]
\item $F(\Phi)$ is a smeared average of the form
$T_f(\Phi) :=   \int_{\mathbb R^d} f(x) \Phi(x) dx
$
 associated with continuous functions of compact support $f$.
 \item $ H(\Phi) $ is a Hamiltonian $H(\Phi):\approx  (\Phi, A  \Phi) + \int_{\mathbb R^d}   P(\Phi(x)  )  \, dx$
with  $(\Phi, A  \Phi)$ a positive definite and reflection-positive (see Section~\ref{sec:6.2}) quadratic form, and 
$P(\Phi(x)  )$ an even polynomial whose terms of  order $\Phi(x)^{2k}$ are interpreted heuristically as representing $k$-particle interactions.  
\end{itemize}
By linearity, the expectation values of products of such variables can be rewritten as
\begin{equation} \label{eq:Schwinger} 
\Big\langle \prod_{j=1}^{n} T_{f_j}(\Phi) \Big\rangle   :=   \int_{(\mathbb R^d)^n}  S_n(x_1,\dots, x_n) \,\prod_{j=1}^n  f(x_j)\,d x_1 \dots d x_n  , 
\end{equation} 
where the $S_n(x_1,\dots, x_n)$ are the {\em Schwinger functions} of the corresponding Euclidean field theory which can be interpreted heuristically as pointwise correlations 
$
\langle \prod_{j=1}^{n} \Phi(x_j) \rangle.
$
Interpreting \eqref{Phi_EV} raises a number of problems of varying difficulty.  

The simplest example of Euclidean fields are the reflection-positive (see Section~\ref{sec:6.2} again) Gaussian fields, for which $H(\Phi)$ contains  only quadratic terms. Gaussian fields are alternatively characterised by  $2n$-point Schwinger functions satisfying Wick's law:
\begin{equation} \label{triv}
S_{2n}(x_1,\dots,x_{2n})  =  \sum_{\pi \text{ pairings}}  \prod_{j=1}^n S_2(x_{\pi(2j-1)},x_{\pi(2j)}).
\end{equation} 
The field theoretical interpretation of \eqref{triv} is the absence of interaction.  Due to that and to their algebraically simple structure, such fields are referred to as {\em trivial}. 

The next level of difficulty is to add the next lowest order even term, i.e.~$\lambda \Phi^4$ for $\lambda>0$.   
Note that, if it exists at all, the corresponding field is a random distribution so making sense of this fourth power is not straightforward. The heuristic RG approach to the problem by Wilson~\cite{Wil71} indicates that in low enough dimensions,  the problem could be tackled through a renormalisation procedure. 
The CQFT  
program  has successfully yielded non-trivial scalar field theories over $\mathbb R^2$ \cite{GliJaf12} and $\mathbb R^3$~\cite{GliJaf73,FelOst76}
, and is still a lively field of mathematical physics.

A natural example aimed at constructing  a  $\Phi^4_d$  functional integral is to  regularise it with a pair of cutoffs: at a short distance (\emph{ultraviolet}) scale and a large distance (\emph{infrared}) scale.    A lattice version of that is the  restriction of  $\Phi(\cdot)$ to  the vertices of a finite graph  
$\Lambda_R^{(a)} := (a \mathbb Z)^d \cap [-R,R]^d$, where $a$ and $R$ play respectively the roles of the ultraviolet and infrared cutoffs.  
For the corresponding finite collection of variables $(\phi_x:x\in \Lambda_R^{(a)})$, the Hamiltonian is then interpreted in terms of a Riemann-sum style discrete analog of the integral expressions, leading to the following statistical-mechanics Gibbs equilibrium state average   
\begin{equation}\label{Phi_SM}
\langle F(\phi) \rangle  = \frac 1 {\mathrm{norm}}   \int_{\mathbb R^{\Lambda_R^{(a)}}}  F(\phi) \exp{[-H(\phi)]}  \prod_{x\in \Lambda_R^{(a)}} d\rho(\phi_x) ,
\end{equation} 
with a Hamiltonian  $H(\phi)$ and an a-priori measure $\rho$   of the form
\begin{equation}\label{H}
H(\phi) = -\sum_{\{x,y\}\subset E(\Lambda_R^{(a)})}\,\phi_x \phi_y    \, , \qquad d\rho(  \phi_x) =  e^{- \lambda \phi_x^4 - b \phi^2_x}  d \phi_x \,,  
\end{equation} 
where $d\phi_x$ is the Lebesgue measure on $\mathbb R$. This is called the {\em $\phi^4$ lattice model}.

The cutoffs are removed through the limit   $R\nearrow \infty$ 
followed by $a\searrow 0$.  Parameters may be added to adjust in the process the spin-spin correlations $\langle\phi_{x_1}\dots\phi_{x_n}\rangle$ in such a way that they stabilise to the Schwinger functions $S_n(x_1,\dots, x_n)$ in the continuum limit scale. 

The Ising model can be thought of as a limiting case of a $\phi^4$ lattice  model as it is obtained by letting $\lambda=b/2$ tend to infinity (the limit of the measures $\rho$ then forces the spins $\phi_x$ to take the values $\pm1$). Actually, the discrete approximations of the $\phi^4$ functional integral and the Gibbs states of an Ising model are {\em always} connected.   This relation is based on a construction which was initiated by Griffiths to obtain the Lee-Yang theorem for the $\phi^4$ lattice  models, and was advanced further by Griffiths and Simon~\cite{GriSim73}.    A probability measure on  $\rho(d\phi) $ on $\mathbb R$ is said to belong to the {\em Griffiths-Simon class} if the expectation values with respect to $\rho$ can be represented as an Ising model on the complete graph with well-chosen coupling constants, or as a limit of such models (satisfying some mild tail conditions). The $\phi^4$ lattice model belongs to the Griffiths-Simon class. For this reason, most techniques that are at our disposal for the Ising model apply to the Griffiths-Simon class. This makes the Ising model an object of major interest when working on CQFT. The developments of the model have therefore been deeply connected to CQFT in the eighties, and we now discuss some examples of such interactions.
 
\subsection{Reflection positivity}\label{sec:6.2}

The notion of reflection positivity was introduced in Quantum Field Theory in the work of Osterwalder-Schrader \cite{OS73}, and we refer to \cite{Bis09} for a review. 
While reflection positivity did not emerge initially as a property of the Ising model, the model remains one of the most natural instances of a reflection positive model, and some of the most striking applications of reflection positivity are indeed dealing with the Ising model.  

Consider the Ising model on a $d$-dimensional torus $\mathbb T_L:=(\mathbb Z/L\mathbb Z)^d$ with $L$ even and split equally the torus into two pieces $\mathbb T_L^+$ and $\mathbb T_L^-$ using hyperplanes (the two pieces are isomorphic to $[0,L/2]\times(\mathbb Z/L\mathbb Z)^{d-1}$) and consider a reflection $\vartheta$ with respect to one of these hyperplanes mapping $\mathbb T_L^+$ to $\mathbb T_L^-$.  We say that $\langle\cdot\rangle$ is {\em reflection positive}  if for all $f,g:\mathbb T_L^+\rightarrow \mathbb R$, 
\begin{equation}
\langle f\vartheta g\rangle=\langle g\vartheta f\rangle \qquad\text{and}\qquad\langle f\vartheta f\rangle\ge 0,
\end{equation}
or, in other words, that $f,g\mapsto \langle f\vartheta g\rangle$ is a positive semi-definite symmetric bilinear form.
The archetypical examples of reflection positive measures are the Ising n.n.f.~measures $\langle\cdot\rangle_{\mathbb T_L,\beta,0}$, but many other examples exist, including some Ising models with long-range interactions.

Reflection positivity has two important implications, namely {\em gaussian domination} leading to the {\em infrared bound}, and the {\em chessboard estimate}. By lack of space, and since most of the applications of reflection positivity to the specific example of the Ising model rely on the infrared bound, let us focus on it and gaussian domination. 

Gaussian domination is a statement linking the partition function of the Ising model with magnetic field to the partition function of the model without it. Formally, it states that  for every function $\mathbf h:V\rightarrow \mathbb R$, $Z_L(\mathbf h)\le Z_L(0)$,
where \begin{equation}
Z_L(\mathbf h):=\sum_{\sigma\in\{-1,1\}^V}\exp\Big[-\beta \sum_{\{x,y\}\in E(\mathbb T_L)}(\sigma_x-\sigma_y+\mathbf h_x-\mathbf h_y)^2\Big].
\end{equation}
Gaussian domination can be proved via reflection positivity through the two hyperplanes mentioned above to show that for each $\mathbf h$, a symmetric version of $\mathbf h$ with respect to a hyperplane has a larger value of $Z_L(\cdot)$. Gaussian domination immediately implies a Fourier version of the infrared bound by using a second-order expansion of $Z_L(\mathbf h)$ near 0: for $d>2$ and every $(a_x)\in \mathbb C^{\mathbb T_L}$ summing to zero,
\begin{equation}
\sum_{x,y\in \mathbb T_L}a_x\overline a_y \langle\sigma_x\sigma_y\rangle_{\mathbb T_L,\beta,0}\le \tfrac2\beta\sum_{x,y\in \mathbb T_L}a_x\overline a_y G(x,y),
\end{equation}
where $G(x,y)$ is the Green function of the simple random walk on $\mathbb Z^d$.

In the specific case of the Ising model, the Messager-Miracle-Sol\'e inequality enables to turn this Fourier estimate into a pointwise estimate on the two-point function: there exist $C,C'>0$ such that for every $\beta>0$ and every $x,y\in \mathbb Z^d$,
\begin{equation}
\langle\sigma_x\sigma_y\rangle_{\beta,0}-m^*(\beta)^2\le \tfrac{C}{\beta}G(x,y)\le \tfrac{C'}{\|x-y\|_2^{d-2}}.
\end{equation} 
This is particularly interesting when $\beta$ approaches $\beta_c$ from below, as it implies that the spin-spin correlations decay algebraically fast at $\beta_c$, with an exponent at least $d-2$.

\subsection{The random current revolution}\label{sec:6.3}

The context of CQFT was also at the origin of one of the most important revolutions in our understanding of the Ising model that we will describe in Section~\ref{sec:6.3.5}. The technique, called the {\em random current representation}, was introduced by Aizenman (inspired by combinatorial identities from \cite{GriHurShe70}). It became one of the most powerful and robust tools available to mathematicians to study the Ising model. We describe it now (see \cite{Dum16} for a review). 

The whole story starts with the observation that the component $\exp[\beta\sigma_x\sigma_y]$ of the Hamiltonian term attached to each edge can be rewritten using Taylor's expansion to get
\begin{equation}
Z(G,\beta,0)=\sum_{\sigma\in\{-1,1\}^V}\prod_{\{x,y\}\in E}\sum_{{\mathbf n}_{\{x,y\}}=0}^\infty \frac{(\beta\sigma_x\sigma_y)^{{\mathbf n}_{\{x,y\}}}}{{\mathbf n}_{\{x,y\}}!}=\sum_{{\mathbf n}\in\mathbb Z_+^E}w_\beta({\mathbf n})\sum_{\sigma\in\{-1,1\}^V}\prod_{x\in V}\sigma_x^{\Delta_x(\mathbf n)},
\end{equation}
where 
\begin{equation}
w_\beta(\mathbf n):=\prod_{\{x,y\}\in E} \frac{\beta^{\mathbf n_{\{x,y\}}}}{\mathbf n_{\{x,y\}}!}\qquad\text{ and }\qquad \Delta_x(\mathbf n):=\sum_{y\in V:\{x,y\}\in E} \mathbf n_{\{x,y\}}.
\end{equation}
Now, the involutions on spin configurations switching the spins at a vertex immediately imply that the sum on $\sigma$ on the right-hand side is either equal to $2^{|V|}$ if $\Delta_x(\mathbf n)$ is even for all $x\in V$, or 0 otherwise (this seems like a very elementary observation, but it bears at the heart of it the $+/-$ symmetry of the space of possible spins). 

Call a function from $E$ to $\mathbb Z_+$ a {\em current}. A {\em source} of the current will be a vertex $x$ with $\Delta_x(\mathbf n)$ odd. The set of sources 
will be denoted by $\partial \mathbf n$.
The previous discussion and the notation lead to the identity
\begin{equation}
Z(G,\beta,0)=2^{|V|}\sum_{\partial\mathbf n=\emptyset}w_\beta(\mathbf n),
\end{equation}
where from now on we omit to specify that we consider currents when using the notation $\mathbf n$.

A current $\mathbf n$ with $\partial\mathbf n=A$ can be interpreted as the occupation time of a collection of paths pairing vertices of $A$ and loops  -- or equivalently the number of times the collection of paths and loops goes through an edge. The decomposition into loops and paths is not unique, nonetheless it remains interesting to interpret currents in terms of them.
 
Proceeding in a similar fashion with the numerator of the spin-spin correlations, we get that 
\begin{equation}
\langle\sigma_A\rangle_{G,\beta,0}=\frac{\sum_{\partial\mathbf n=A}w_\beta(\mathbf n)}{\sum_{\partial\mathbf n=\emptyset}w_\beta(\mathbf n)}.
\end{equation}
In words, one may write spin-spin correlations in terms of weighted sums of currents with specific source constraints $\partial\mathbf n=A$ and $\partial\mathbf n=\emptyset$.
Note that the source constraint is not the same for the numerator and denominator.

\begin{framed}
\centerline{
\textbf{Frame~4: the high-temperature expansion and $\beta_c>0$}
}
The {\em high-temperature expansion} of the Ising model, due to van der Waerden \cite{Wae41}, can be neatly defined here as the set of edges with an {\em odd} current (it can also be obtained by a direct expansion using that $\exp[\beta \sigma_x\sigma_y]=\cosh(\beta)+\sinh(\beta)\sigma_x\sigma_y$). One ends up with another expression of the partition function in terms of even subgraphs
\begin{equation}
Z(G,\beta,0)=\cosh(\beta)^{|E|}\sum_{F\in \mathrm{Even}(G)}\tanh(\beta)^{|F|},
\end{equation} 
which resembles the low-temperature expansion, except that it is on $G$ instead of $G^*$  and that it is valid for arbitrary graphs and not only planar ones. In particular, one may easily deduce the Kramers-Wannier duality between the low and high temperature expansions at temperatures $\beta$ and $\beta^*$ satisfying $\tanh(\beta)=e^{-2\beta^*}$ in the case of the square lattice.

One application of currents (or alternatively high-temperature expansion) is obtained by considering a mapping from currents with $\partial\mathbf n=\{x,y\}$ to currents with $\partial\mathbf n=\emptyset$ setting the current on a path from $x$ to $y$ of odd current  (such a path necessarily exists) to $0$. This many-to-one mapping (one has to keep track of the path and the value of the current on it to reconstruct the preimage) increases drastically the weight of the current as soon as $\beta\ll1$, which shows that the spin-spin correlations $\langle\sigma_x\sigma_y\rangle_{G,\beta,0}$ are decaying exponentially fast in this regime. This implies in particular that $\beta_c>0$. 
\end{framed}

A key observation of Aizenman is that the so-called {\em switching lemma}, see Frame~5, pertaining to combinatorial properties of the random current model, could be used to reinterpret spin-spin correlations as well as many other properties in terms of {\em probabilities} involving multiple independent currents. This lemma completely changed the point of view on currents, as it transforms them from a combinatorial type object into a probabilistic one. In particular, intuitions coming from probabilistic models such as random walks and percolation was later used to prove new theorems on the Ising model; see Sections~\ref{sec:6.3.5}, \ref{sec:6.5}, and \ref{sec:7.1}.

\begin{framed}
\centerline{
\textbf{Frame~5: the switching lemma for random currents}
}
Write ${\mathbf n}\in \mathcal F_A$ if there exists $\mathbf k\le\mathbf n$ with $\partial \mathbf k=A$. Note that if $A=\{x,y\}$, this is equivalent to the existence of a path from $x$ to $y$ which is made of edges with a positive current. Recall that $A\Delta B$ denotes the symmetric difference of the sets $A$ and $B$. With this notation, the {\em switching lemma} \cite{GriHurShe70} states that for every $F:\mathbb Z_+^E\rightarrow \mathbb R$ and every two sets of vertices $A,B\subset V$, \begin{equation}
\sum_{\substack{\partial\mathbf n_1=A\\ \partial\mathbf n_2=B}}w(\mathbf n_1)w(\mathbf n_2)F(\mathbf n_1+\mathbf n_2)=\sum_{\substack{\partial\mathbf n_1=A\Delta B\\ \partial\mathbf n_2=\emptyset}}w(\mathbf n_1)w(\mathbf n_2)F(\mathbf n_1+\mathbf n_2)\mathbb I(\mathbf n_1+\mathbf n_2\in \mathcal F_B).
\end{equation}
The name of the lemma is fairly self-explanatory, as it consists, when considering sums of two currents, of a recipe to switch the sources from the second one to the first one. The proof is a very entertaining combinatorial problem that is left to the reader. 

A direct application (to illustrate the strength of the lemma) is the case $A=B$, which gives immediately that 
\begin{equation}
\langle\sigma_A\rangle_{G,\beta,0}^2=\mathbb P^\emptyset_G\otimes\mathbb P^\emptyset_G[\mathbf n_1+\mathbf n_2\in \mathcal F_A],
\end{equation}
where $\mathbb P^B_G$ is the measure on currents $\mathbf n$ on $G$ with $\partial \mathbf n=B$ attributing to each such $\mathbf n$ a probability that is proportional to $w(\mathbf n)$, and $\otimes$ denotes the product for probability measures. In words, one may interpret the {\em square of} spin-spin correlations $\langle\sigma_A\rangle_{G,\beta,0}$ as the probability, for the sum of two {\em independent} random currents, of pairing the elements of $A$ by paths of positive current. One may also try as an exercise to recover Griffiths' inequalities from the switching lemma.
\end{framed}

\subsection{Triviality in dimension $d>4$}\label{sec:6.3.5}

In 1982, Michael Aizenman and Juerg Fr\"ohlich \cite{Aiz82,Fro82} independently proved that the scaling limit of the Ising model is trivial in dimension five and more in the following sense. 
Consider discrete smeared averages defined by 
\begin{equation}\label{def_Tf_scaled}
T_{f,L}(\sigma) :=  \frac{1}{\sqrt{\Sigma_L}}\sum_{x\in \mathbb Z^d} f(x/L ) \, \sigma_x   \,,
\end{equation}
where $f$ ranges over compactly supported continuous functions,  and $\Sigma_L:=\big\langle\big(\sum_{x\in\Lambda_L}\sigma_x\big)^2\big\rangle$ denotes the variance of the sum of spins over the box of size $L$.
The theorem states that when $d>4$, these smeared averages $T_{f,L}(\sigma)$ are approximately Gaussian of variance $\langle T_{f,L}(\sigma)^2\rangle_\beta$
in the sense that there exists an explicit constant $C_f>0$ such that for every $\beta\le\beta_c$, every $L\le\xi(\beta)$, and every $z>0$,
\begin{equation}\label{eq:gauss}
 \Big|\,\big\langle \exp[z\,T_{f,L}(\sigma)-\tfrac{z^2}2\langle T_{f,L}(\sigma)^2\rangle_\beta ]\big\rangle_\beta\,-\,1\,\Big|~\le~ \frac{C_fz^4}{L^{d-4}}.
\end{equation}
In words, the previous statement claims that the characteristic function of $T_{f,L}(\sigma)$ is close to the one of a Gaussian random variables.

As a direct consequence of this result, one obtains that any well-defined scaling limit of the Ising model, and in fact more generally of the $\phi^4$ lattice model, is inevitably Gaussian. 
The result marked a brutal stop in the CQFT program outlined in Section~\ref{sec:6.1} as the proofs suggested, while not proving, that the model should also be trivial in four dimensions.

As mentioned above, one of the most striking applications of the random current representation is related to CQFT. Indeed, Aizenman's proof of this theorem relies on a beautiful parallel between random walks and the paths joining sources in currents. We do not resist discussing this link below. But before doing so, let us mention that the approach of Fr\"ohlich in \cite{Fro82}, based on the Brydges-Fr\"ohlich-Spencer (BFS) walk representation of spin-spin correlations \cite{BryFroSpe82}, is deeply connected to the random current as well. The walks in the BFS representation play the roles of the paths between sources in the random current. The advantage of this alternative approach is that it works for more general models, at the cost of losing the switching lemma and its benefits.

Let us focus on the four-point function and define the corresponding Ursell function given, for $x_1,\dots,x_4\in \mathbb Z^d$, by 
\begin{equation}
U_{4}^\beta(x_1,\dots,x_4):=\langle\sigma_{x_1}\cdots\sigma_{x_4}\rangle_\beta-\sum_{\pi\text{ pairing}}\prod_{i=1}^2\langle\sigma_{x_{\pi(2i-1)}}\sigma_{x_{\pi(2i)}}\rangle_\beta.
\end{equation}
A simple exercise involving the switching lemma shows that 
\begin{equation}
U_{4}^\beta(x_1,\dots,x_4)=-2\langle\sigma_{x_1}\sigma_{x_2}\rangle\langle\sigma_{x_3}\sigma_{x_4}\rangle\mathbb P^{\{x_1,x_2\}}\otimes\mathbb P^{\{x_3,x_4\}}[x_1,\dots,x_4\text{ all connected in }\mathbf n_1+\mathbf n_2],
\end{equation}
where connected in $\mathbf n_1+\mathbf n_2$ means being connected by a path of edges with $\mathbf n_1+\mathbf n_2$ not equal to zero. If one remembers that one can think of a current with sources $x_1$ and $x_2$ as a path connecting the two vertices together with a collection of loops, one can reinterpret the right-hand side of the previous identity at the light of so-called random walks (a random walker traces his way through the vertices of a graph by picking its next steps at random among neighbours of where it currently stands -- this Markov process is one of the most fundamental objects of probability theory). It is a classical result that two random walks connecting two pairs of points that are at a mutual distance of order $L$  intersect with a probability bounded away from 0 as $L$ tends to infinity in dimensions $d<4$, and tending to zero in dimension $d\ge4$. 

At this stage, it is totally unclear why the paths linking the points $x_1$ and $x_2$ in $\mathbf n_1$, and $x_3$ and $x_4$ in $\mathbf n_2$, would behave as random walks. It is also unclear what would be the impact  of the additional loops. Still, it is tempting to think that if an analogy with random walks was valid, then it would single out dimensions $d\ge4$ as being dimensions for which $U_4^\beta$ becomes much smaller than products of two-point correlations or, in other words, for which Wick's law would become asymptotically valid, thus hinting at triviality. 

When the dimension is strictly larger than 4, the story for random walks becomes even simpler, as the {\em expected number of intersections} is also tending to zero with $L$. Using the infrared bound to estimate the spin-spin correlations of the Ising model, one may go around the difficulty of proving a random walk type behaviour for currents to show that the intersection probability is tending to 0.

Making the argument work for currents in dimension 4 is more subtle because, contrarily to larger dimensions, 
the expected number of intersections does not tend to 0 when $L$ tends to infinity. Hence, in order to prove that the intersection probability goes to 0, 
one inevitably has to go deeper in the understanding of the analogy between currents and random walks.

\subsection{Rigorous renormalisation group in 4D Ising}\label{sec:6.4}

The triviality of the Ising model in dimension $d>4$ naturally raises the question of its triviality in dimension $d=4$, which is not only the pertinent physical dimension for CQFT, but also for the so-called $4-\varepsilon$ expansions providing information on dimension 3. 
In the eighties, Wilson's renormalisation group method was already in every physicists' toolbox, yet the challenges to overcome to cast the general theory in a mathematical framework seemed out of reach. 
Interestingly, a very relevant case became an important exception. 

Consider the lattice version of the $\phi^4$ model discussed in Section~\ref{sec:6.1}. The case $b=\lambda=0$ corresponds to a Gaussian field known under the name of discrete Gaussian Free Field (GFF), which enjoys a number of striking features. One of them is that the model converges, when rescaling the lattice, to the continuum GFF. In a series of impressive papers \cite{GawKup85,FMRS87, HarTas87}, mathematical physicists proved in the eighties that, when starting from a weakly coupled $\phi^4$ lattice model (meaning that $\lambda$ is small), one may apply a multi-scale analysis to prove convergence of the model to the continuum GFF. 

Several methods were used at the time, but let us mention that the method of Gawedski and Kupiainen \cite{GawKup85} can be thought of as a rigorous version of Kadanoff block-spin renormalisation procedure. It consists of writing the model in terms of averages of spins over large blocks of size $L^k$, and to average them out scale by scale. At leading order, each step of the procedure boils down to modifying the parameters of the model. Of course, the reality is much more complicated than the first order analysis suggests, and the renormalisation scheme is quite complex.

An alternative to this block-spin renormalisation was later developed by Bauerschmidt, Brydges, and Slade \cite{BauBrySla14} in order to obtain refined results, as well as to treat more general models. In these alternative approaches, the block-spin analysis is replaced by the following strategy: one thinks of quantities in the $\phi^4$ lattice model as being expressed in terms of the discrete GFF itself. In order to control the asymptotic behaviour of such quantities, one decomposes the covariance of the discrete GFF into a sum of finite-range covariances that one integrates out one by one. At each step a change of the parameters of the system is required to keep things converging towards a limit. Doing so enables the authors to focus their attention on how the parameters evolve under this procedure. This evolution can be thought of as the renormalisation map in the renormalisation group. 

The level of sophistication of these techniques is quite astonishing, and the precision of the results outstanding. As one may guess, this comes at a price. At the bottom of both strategies lies the fact that the original $\phi^4$ lattice model is in the ``vicinity of a model'', the Gaussian Free Field, that enjoys a number of nice properties. As a result, the technique is (as for today) {\em perturbative} in nature, which is somehow its main limitation. We will see another instance of such a renormalisation scheme, this time near another fixed point, when discussing the 2D Ising model.

\subsection{Forty years later: the random current strikes back}\label{sec:6.5}

While renormalisation techniques provided impressive rigorous results in dimension 4, they remained as we mentioned perturbative, meaning that they required that the lattice $\phi^4$ model one starts from has a small $\phi^4$ term. Yet, if one would like to construct a non-trivial 4D quantum field theory, one would definitely try to start with a strongly coupled $\phi^4$ lattice model (meaning with a $\phi^4$ terms which is not a priori small), for instance working with the Ising model which in some sense can be thought of as the model with the strongest possible coupling, thus excluding existing renormalisation group techniques. 

This asks for another approach, and this is probably why one had to wait for forty years to finally obtain a proof of the triviality of the 4D Ising and $\phi^4$ lattice models, which states \cite{AizDum21} that there exists $c>0$ such that for the n.n.f.~$\phi^4$ lattice model on $\mathbb Z^4$ with parameters $b,\lambda$, and a compactly supported continuous function $f$, there exists $C_f>0$ such that for every $\beta\le\beta_c=\beta_c(b,\lambda)$, every $L\le\xi(\beta)$, and every $z>0$,
\begin{equation}\label{eq:gauss}
 \Big|\,\big\langle \exp[z\,T_{f,L}(\varphi)-\tfrac{z^2}2\langle T_{f,L}(\varphi)^2\rangle_\beta ]\big\rangle_\beta\,-\,1\,\Big|~\le~ \frac{C_fz^4}{(\log L)^c}.
\end{equation}

The strategy of the proof uses a more delicate probabilistic perspective on the random current than in \cite{Aiz82}, still keeping in mind the interpretation in terms of random walks of the paths joining the sources of the current. Indeed, it can be proved that two random walkers in
four dimensions going from points to points that are all at a mutual distance of order $L$ intersect with probability of order $(\log L)^{-c}$ for some universal constant $c>0$. The reason is that while the expected number of intersections is of order 1, the number of intersections, when such intersections exist, is with high probability quite large in $L$ (and is growing with $L$). The core of the paper is to apply a similar argument to the paths in the random current. Of course, challenges emerge when trying to handle the highly non-Markovian paths obtained by considering the paths joining the sources in currents. Nevertheless, guided by the random walk intuition, one can build a multi-scale analysis to prove that conditioned on intersecting, random currents intersect a large number of times, and ultimately deduce from this the triviality result.

\section{The last fifty years: Ising model and percolation}\label{sec:7}

Percolation theory gathers under its umbrella a variety of random graphs systems.  A configuration on $G=(V,E)$ is an element $\omega=(\omega_e:e\in E)\in\{0,1\}^E$ which is interpreted as a subgraph  with vertex-set $V$ and edge-set $\{e\in E:\omega_e=1\}$. Then, different percolation  models can be defined by considering different measures on $\{0,1\}^E$. Historically, the original model, called Bernoulli percolation, is defined in such a way that the $\omega_e$ are independent Bernoulli random variables. It was introduced to understand the behaviour of liquid in a porous medium. Nevertheless, the theory of non-Bernoulli models has been found to be related to a variety of other models of statistical physics explaining various physical phenomena. 

As often, the Ising model has played an essential role in the development of percolation theory, and conversely certain advances in percolation theory have been fundamental to our understanding of the Ising model.
Sometimes, the link between the two models is simply an analogy between their behaviours, but sometimes the connection is much more direct. For instance, spin-spin correlations can be rewritten in terms of a percolation model, in which case we speak of the percolation model as being a {\em graphical representation} of the Ising model.  We now propose to discuss some examples of these links between the Ising model and percolation.

\subsection{Percolation interpretation of random currents}\label{sec:7.1}

We have seen one example of a graphical representation in Frame~5 where the squares of spin-spin correlations get rephrased as connectivity properties of the sum of two currents. One may easily define a percolation model out of the pair of currents above by saying that for an edge $\{x,y\}$, $\omega_{\{x,y\}}=1$ if $(\mathbf n_1+\mathbf n_2)_{\{x,y\}}>0$. Then, the square of the spin-spin correlations between two points becomes the probability, for this percolation model, that $x$ and $y$ are connected in $\omega$.  

The best illustration of how intuition from percolation or the Ising model can drive developments on the other model is provided by an important result on the Ising model in the regime $\beta<\beta_c$. This result from 1987, due to Aizenman, Barsky and Fernandez \cite{AizBarFer87} (see \cite{DumTas15} for an alternative argument), states that correlations of the n.n.f.~Ising model decay exponentially fast as soon as $\beta<\beta_c$ in the sense that for each such $\beta$, there exists $\tau>0$ such that for every $x,y\in \mathbb Z^d$,
\begin{equation}\label{eq:exponential decay}
\langle\sigma_x\sigma_y\rangle_{\beta,0}\le \exp(-\tau\|x-y\|).
\end{equation}
We say that the phase transition is {\em sharp}: there is no intermediate phase $(\beta_{\mathrm{exp}},\beta_c)$ in the Ising model in which spin-spin correlations would decay polynomially. Let us mention that a similar exponential decay was obtained recently for truncated correlations $\langle\sigma_x\sigma_y\rangle_{\beta,0}-m^*(\beta)^2$ when $\beta>\beta_c$, see \cite{DumGosRao18}.

This theorem is of fundamental importance for the following reason. Perturbative results, which are combinatorial in nature, are valid under the assumption that certain quantities decay exponentially fast, and in fact with a rate of decay which is sufficiently large. While this hypothesis is important to apply the techniques, it happens to be of little relevance from a physical point of view. In fact, one expects that most of the phenomenology remains unchanged as long as spin-spin correlations decay exponentially fast. As a consequence, \eqref{eq:exponential decay} can be thought of as a bottleneck in the understanding of the phase $\beta<\beta_c$: as soon as it is obtained, a number of important results can be derived from it. As an example, the results on fluctuations of interfaces and Ornstein-Zernike estimates were proved to hold in the whole regime $\beta<\beta_c$. The result also provides meaning to the correlation length $\xi(\beta)$ mentioned in Section~\ref{sec:4.1}, as it proves that it is finite as soon as $\beta<\beta_c$.

Let us now comment on the proof. The argument relies on a fruitful idea consisting in deriving differential inequalities between thermodynamical quantities of the Ising model. The archetypical example of such differential inequalities are given, for the problem at hand, by (recall that the magnetisation $m=m(\beta,h)$ is a function of $\beta$ and $h$)
\begin{equation}
m\le \tanh(\beta h)\tfrac{\partial }{\partial(\beta h)}m+m^2(\beta\tfrac{\partial }{\partial\beta}m+m)\qquad\text{ and }\qquad m\tfrac{\partial }{\partial\beta}m\ge c. 
\end{equation}
The interesting feature here is that similar differential inequalities appear when studying Bernoulli percolation. In fact, a number of results were obtained in parallel during the eighties, where each result for Ising had its pendant for Bernoulli percolation, and vice versa. As an example, critical exponents for $d>4$ were obtained by Aizenman and Fernandez \cite{AizFer86} using differential inequalities that can be adapted to Bernoulli percolation. These techniques are useful to transform qualitative results (e.g.~a quantity tends to 0) to quantitative ones (e.g.~exponentially fast). We do not resist mentioning one of them: for $h=0$ and $\beta<\beta_c$,
\begin{equation}
\Big(1-\frac{B}{\chi}\Big)\frac{2d\chi^2}{1+B}\le\tfrac{\partial}{\partial\beta}\chi\le 2d\chi^2,
\end{equation}
where $\chi(\beta):=\sum_x\langle\sigma_0\sigma_x\rangle_{\beta,0}$ is the {\em susceptibility}, and $B(\beta)$ is the {\em Bubble diagram} and is given by
\begin{equation}
B(\beta):=\sum_{x\in \mathbb Z^d}\langle\sigma_0\sigma_x\rangle_{\beta,0}^2.
\end{equation}
Since the Infrared Bound implies that $B(\beta)$ remains bounded uniformly in $\beta<\beta_c$ as soon as $d>4$, $\chi(\beta)$ must blow up like $1/|\beta-\beta_c|$ as $\beta$ approaches $\beta_c$ from below.

Another striking instance of how fruitful the connection between percolation models and the Ising model was for the development of both models is the following {\em continuity result} of the phase transition of the 3D Ising model, due to \cite{AizDumSid15}, stating that  the n.n.f.~Ising model satisfies $m^*(\beta_c)=0$ for every $d\ge3$.

The argument relies on percolation methods applied to the double random current representation of an argument of Burton and Keane proving the uniqueness of the infinite connected component of percolation. The whole argument can be improved and extended to study all translation-invariant Gibbs measures, obtaining the classification result already mentioned in Section~\ref{sec:5.2}. 

\subsection{Fortuin-Kasteleyn percolation}\label{sec:7.2}

  Another (and in fact older) example of a graphical representation is provided by a special case of the Fortuin-Kasteleyn (FK) percolation. In this model, introduced in \cite{ForKas72}, the measure $\phi_{G,p,q}$ is given, for $G=(V,E)$ finite and  $\omega\in\{0,1\}^E$, by 
\begin{equation}
\phi_{G,p,q}[\{\omega\}]:=\frac1{Z(G,p,q)}p^{|\omega|}(1-p)^{|E|-|\omega|}q^{k(\omega)},
\end{equation}
where $p\in[0,1]$ and $q>0$ are the parameters of the model, called respectively the {\em edge-weight} and the {\em cluster-weight}, $|\omega|:=\sum_{e\in E}\omega_e$ is interpreted as the number of edges in $\omega$, and $k(\omega)$ is the number of connected components of $\omega$. 

When $q=1$, one ends up with the classical Bernoulli percolation model in which the $\omega_e$ are independent. When $q\ne1$, the state of edges is no longer independent and one ends up with a {\em dependent percolation model} whose study is central in modern probability theory. From now on, we focus on the case $q=2$, which we call the {\em FK Ising model}.
We confine our discussion to two features of this percolation model, namely its link to the Ising model, and the FKG inequality.

\begin{figure}[t]
\begin{center} \includegraphics[width=1.0\textwidth]{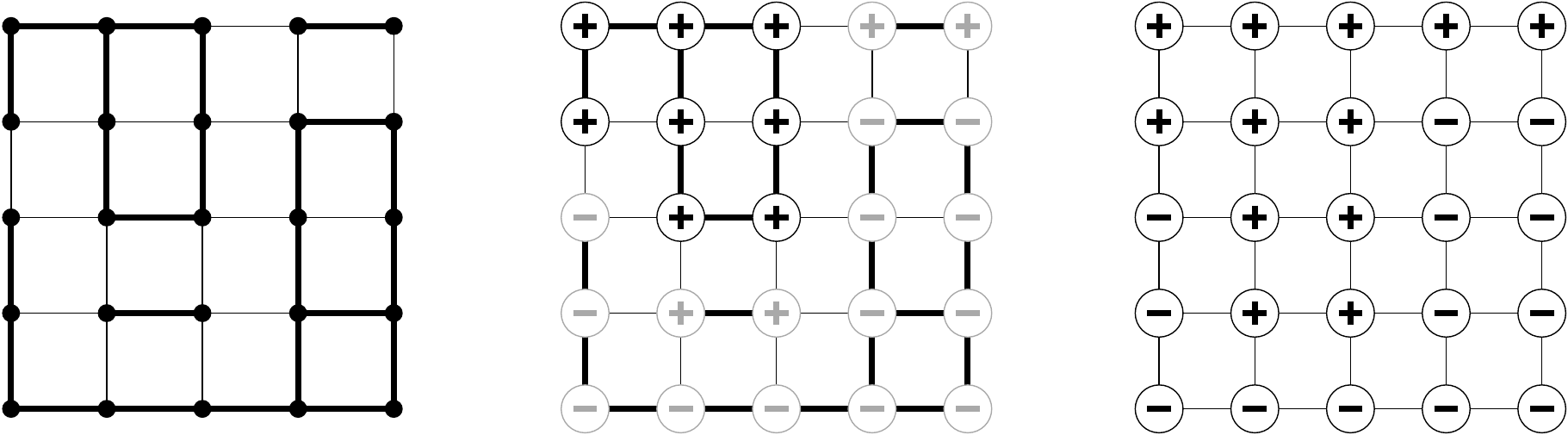} 
\caption{The Edwards-Sokal coupling, with on the left a picture of the FK Ising configuration (bold edges are those with $\omega_e=1$), in the middle, spins are attached to each cluster (one example in black and others in grey), and on the right, the spins without the FK Ising configuration.}\label{fig:3}\end{center}
\end{figure}

Let us start with the former, which provides a recipe to obtain the Ising model configuration out of FK Ising; see Figure~\ref{fig:3}. 
Consider a random variable $\omega\in\{0,1\}^E$ with the law of FK Ising with parameter $p\in[0,1]$ and construct $\sigma\in\{-1,1\}^V$ by 
\begin{itemize}[noitemsep,nolistsep]
\item choosing for every connected component $\mathcal C$ of $\omega$ a spin $\sigma_\mathcal C$ uniformly between $-1$ and $+1$, and independently of the other connected components. 
\item defining $\sigma_x=\sigma_\mathcal C$ for every $\mathcal C$ and every $x\in \mathcal C$.
\end{itemize}
Then, $\sigma$ has the law of the Ising model on $G$ with parameter $\beta=\tfrac12\log[1/(1-p)]$ and $h=0$. This coupling, due to Fortuin and Kasteleyn and often referred to as the {\em Edwards-Sokal coupling} due to the paper \cite{EdwSok88},  enables to express correlation functions of the Ising model in terms of FK Ising. For instance, by decomposing on the events that $x$ is connected to $y$ or not in $\omega$, one easily gets that
\begin{equation}
\langle\sigma_x\sigma_y\rangle_{G,\beta,0}=\phi_{G,1-e^{-2\beta},2}[x\text{ connected to }y\text{ in }\omega].
\end{equation}
Similarly, $\langle\sigma_A\rangle_{G,\beta,0}=\phi_{G,1-e^{-2\beta},2}[\mathcal F_A]$, where $\mathcal F_A$ is the event that each connected component of $\omega$ contains an even number (possibly equal to 0) of vertices in $A$. Another interesting feature of this coupling is  that it is at the basis of so-called {\em cluster algorithms} due to Swendsen and Wang,  who used it to speed up the Glauber dynamics and the simulation of the Ising model, in particular near the critical point.
The interest of FK Ising and more generally FK percolation models {\em with $q\ge1$} is that they enjoy some nice monotonicity properties (dependent percolation models satisfying these properties have been an object of intense study in the past ten years). 
 Let us mention two such properties.  The Fortuin-Kasteleyn-Ginibre (FKG) inequality states that for every increasing functions $f,g:\{0,1\}^E\rightarrow \mathbb R$,
\begin{equation}
\phi_{G,p,q}[fg]\ge \phi_{G,p,q}[f]\phi_{G,p,q}[g].
\end{equation}
This inequality is often used for indicator functions of {\em increasing events} (i.e.~events for which the indicator function is an increasing function), in which case the inequality states that increasing events are positively correlated. Another manifestation of the monotonicity properties is the monotonicity in $p$: for every increasing function $f:\{0,1\}^E\rightarrow \mathbb R$ and $p'\ge p$,
\begin{equation}
\phi_{G,p',q}[f]\ge \phi_{G,p,q}[f].
\end{equation} 

These monotonicity properties are particularly useful. The second one applied to FK Ising and the indicator function of $\mathcal F_A$ implies that $\langle\sigma_A\rangle_{G,\beta,0}$ is increasing in $\beta$, and the first one applied to indicator functions of $\mathcal F_A$ and $\mathcal F_B$ implies the second Griffiths inequality $\langle\sigma_A\sigma_B\rangle_{G,\beta,0}\ge \langle\sigma_A\rangle_{G,\beta,0}\langle \sigma_B\rangle_{G,\beta,0}$.

\subsection{The broader impact of the Ising model on dependent percolation models}\label{sec:7.3}

In the first fifty years that followed its introduction, the theory of percolation was much more advanced for Bernoulli percolation than for other dependent percolation models. The past ten years have seen tremendous progress in bridging the gap between our understanding of the Bernoulli case and the others. The interplay between dependent percolation models and the Ising model has been fundamental for these developments. 

We already saw that the Ising model is related to FK Ising and a percolation model created out of random currents. It does not come as a surprise that one of the first dependent percolation models to see significant progress in its understanding was the FK Ising. Of course, the Edwards-Sokal coupling enables to transfer immediately certain known facts about the Ising model to its percolation representation (for instance, the critical point of the FK Ising on $\mathbb Z^2$ is $1-e^{-2\beta_c}=\sqrt2/(1+\sqrt 2)$ thanks to Onsager's result). Also, the model enjoys some specific features that make its direct analysis simpler than for other dependent percolation models. 

For all these reasons, the FK Ising became the entrance gate to a new realm of results on dependent percolation models. A perfect illustration of this is provided by the study of {\em crossing probabilities} for planar dependent percolation models. Let us provide slightly more detail. 

One important feature of critical dependent percolation models in two dimensions is that they satisfy the box-crossing property (BCP), and its connected notion the  Russo-Seymour-Welsh theory (RSW). More precisely, if for a rectangle $R$, the event $\mathrm{Cross}(R)$ corresponds to the existence of a path in $\omega$ between the left and right sides of $R$, the properties (BCP) and (RSW) for a percolation model on $\mathbb Z^2$ with measure $\mathbb P$ are the following:
 \begin{itemize}
 \item (BCP) for all $\rho>0$, there exists $c>0$ such that for every $n\ge1$,
 \begin{equation}
 c\le \mathbb P[\mathrm{Cross}([0,\rho n]\times[0,n])|\omega_{|[-n,(\rho+1)n]\times[-n,2n]^c}]\le 1-c \text{ almost surely}.\end{equation}
 \item (RSW) for all $\rho>0$, there exists $C>0$ such that for every $n\ge1$, 
 \begin{equation}
 \mathbb P[\mathrm{Cross}([0,\rho n]\times[0,n])]\ge \mathbb P[\mathrm{Cross}([0,n]\times[0,\rho n])]^C.
\end{equation}
\end{itemize}
These two properties have been the driving force of the progress in our understanding of the 2D dependent percolation models. The FK-Ising model played an essential role in these developments, as it was the first dependent percolation model for which (BCP) could be proved \cite{DumHonNol11}. This development triggered a whole new direction of research that led to substantial progress in our understanding of (BCP) and (RSW) for various percolation models.

\section{Over the last ten years: conformal invariance of the Ising model}\label{sec:8}

\subsection{What is conformal invariance?}\label{sec:8.1}

As mentioned before, Kadanoff used his block-spin renormalisation to predict that the large scale properties of the critical Ising model were invariant under scaling. The same argument also leads to postulate {\em translation} and {\em rotation invariance}. In 1970, Polyakov \cite{Pol70}
suggested a much stronger invariance of the model. Since we saw that it is natural to associate a QFT with the large scale properties of the critical Ising model, and since this QFT is a local field, these properties should be invariant under any map which is {\em locally} a composition of translation, rotation and homothety. As a corollary one predicts full {\em conformal invariance}, i.e.~invariance  under all one-to-one holomorphic maps. This prediction was turned into a classification of possible {\em conformal field theories} (CFT) in 2D in seminal papers by Belavin, Polyakov and Zamolodchikov \cite{BelPolZam84} that generated an explosion of activity, allowing non-rigorous explanations of many critical phenomena.

From a mathematical perspective, the notion of conformal invariance of a model is  not straightforward to define. A number of interpretations of the limit of large scale properties -- called the {\em scaling limit} -- can be taken, and we mention a few now. 

For clarity of the exposition, we focus on the critical Ising model on $\mathbb Z^d$ and its rescaled versions $a\mathbb Z^d$ for $a>0$. We drop the subscript referring to $\beta$ and $h$ as they are fixed to be equal to $\beta_c$ and $0$ respectively. Consider a simply connected domain $\Omega\subsetneq \mathbb R^d$. 
\medbreak\noindent
\textbf{(Spins)} The most natural approach is to consider the spin-spin correlations  defined for every $a>0$ and $x_1,\dots,x_n\in \Omega$ by
\begin{equation}
S_{\Omega}^{(a)}(x_1,\dots,x_n):=\langle\sigma_{[x_1]_a}\dots\sigma_{[x_n]_a}\rangle_{a\mathbb Z^d\cap\Omega},
\end{equation} where $[x]_a$ is the vertex of $a\mathbb Z^d\cap\Omega$ closest to $x$. These Schwinger functions already appeared as the key players in CQFT. 
One is then interested in the limit as $a$ tends to 0 of these properly renormalised quantities. If the limit exists, we call it $S_\Omega(x_1,\dots,x_n)$. 
\medbreak\noindent
\textbf{(Energies)} Another object of interest is the {\em energy-energy} correlations. For $a>0$ and $x_1,\dots,x_n\in \Omega$, one considers at the quantities
 \begin{equation}
T_{\Omega}^{(a)}(x_1,\dots,x_n):=\langle\varepsilon_{(x_1)_a}\dots\varepsilon_{(x_n)_a}\rangle_{a\mathbb Z^d\cap\Omega},
\end{equation} where $\varepsilon_{\{u,v\}}:=\sigma_u\sigma_v-\langle\sigma_u\sigma_v\rangle_{a\mathbb Z^d}$ and $(x)_a$ is the edge closest to $x$. The quantity $\varepsilon_x$ is called the {\em energy}. One is again interested in the limit $T_\Omega(x_1,\dots,x_n)$ as $a$ tends to 0 of these properly rescaled quantities. 
\medbreak\noindent
\textbf{(Geometry of interfaces)} In two dimensions, another direction was proposed in the nineties. It consists in considering the low-temperature representation, i.e.~the interfaces between plus and minus spins. In a domain $\Omega$, it creates a family of non-intersecting loops together with arcs from boundary to boundary. Let $\mathfrak C_\Omega$ be the set of such collections of loops and arcs. The set $\mathfrak C_\Omega$ can be turned into a metric space by attaching a distance $d_\Omega$ which, heuristically, states that two configurations are close to each other when the large loops and arcs are close to each other. Let us call $\mathcal C^{(a)}_\Omega$ the random variable obtained by considering the low-temperature expansion of a critical Ising model configuration in $a\mathbb Z^d\cap \Omega$. Here, we are interested in the limit of $\mathcal C^{(a)}_\Omega$ as a random object.

Now, what do we mean by conformal invariance? Roughly speaking, we mean that certain quantities of the model are conformally covariant/invariant. With the definitions above, it would for instance mean that there exists a way of renormalising the $S_\Omega^{(a)}(x_1,\dots,x_n)$ and $T_\Omega^{(a)}(x_1,\dots,x_n)$ in such a way that they converge to quantities $S_\Omega(x_1,\dots,x_n)$ and $T_\Omega(x_1,\dots,x_n)$ that satisfy that there exist $\bm\Delta_\sigma,\bm\Delta_\varepsilon$ such that for every conformal (i.e.~holomorphic and one-to-one) map $f:\Omega\rightarrow f(\Omega)$, we have
\begin{align}
S_{f(\Omega)}(f(x_1),\dots,f(x_n))&=|f'(x_1)|^{-\bm\Delta_\sigma}\cdots|f'(x_n)|^{-\bm\Delta_\sigma}\, S_\Omega(x_1,\dots,x_n),\\
T_{f(\Omega)}(f(x_1),\dots,f(x_n))&=|f'(x_1)|^{-\bm\Delta_\varepsilon}\cdots|f'(x_n)|^{-\bm\Delta_\varepsilon}\, T_\Omega(x_1,\dots,x_n).
\end{align}
For the geometry of interfaces, the situation is even simpler as one means that the family of loops and arcs $\mathcal C_\Omega^{(a)}$ converges to a limit $\mathcal C_\Omega$ as $a$ tends to $0$ and that this limit satisfies that 
$\mathcal C_{f(\Omega)}$ and $f(\mathcal C_\Omega)$ have the same law for every conformal map $f:\Omega\rightarrow f(\Omega)$.

\subsection{Conformal invariance of the 2D Ising model}\label{sec:8.2}

Around fifteen years ago, Smirnov \cite{Smi10} and Chelkak and Smirnov \cite{CheSmi12} obtained a major breakthrough towards proving conformal invariance of 2D Ising model. This fundamental proof, that we discuss below, opened the way to a very deep understanding of the scaling limit of the model. 

A few years later, Chelkak-Izyurov-Hongler \cite{CheHonIzy15} proved conformal covariance of the spin-spin correlations (with $\bm\Delta_\sigma=1/8$). It was later proved in \cite{CamGarNew15} that the quantities $S_\Omega(x_1,\dots,x_n)$ are the Schwinger functions of a random distribution, that can be understood as the {\em spin-field} that physicists sometimes refer to. In the same spirit, conformal covariance of the energy-energy correlations was proved in \cite{Hon10} (with $\bm\Delta_\varepsilon=1$). In this case, one may prove that the correlations are {\em not} the Schwinger functions of a  random distribution. Turning to interfaces, the following result was the culmination of the theory: the arcs in $\mathcal C_\Omega$ are given by the so-called free arc ensemble of parameter 3 and the loops by conformal loop ensembles of parameter $3$ in the simply connected domains obtained as the complements of the arcs (see \cite{BenDumHon14,BenHon16}). In particular, the scaling limit is conformally invariant. This body of work uses the ideas from \cite{Smi10,CheSmi12} together with the theory of the Schramm-Loewner evolution and its consequences.

As mentioned above, an important breakthrough came from the works \cite{Smi10,CheSmi12} where conformal covariance of so-called {\em fermionic observables} $f_\Omega^{(a)}$ is proved. Those observables are linear combinations of order-disorder operators  (see Frame~6) considered by Kadanoff and Ceva in \cite{KadCev71}, see also \cite{CheCimKas15} for several connections to other classical objects.

\begin{framed}
\centerline{\textbf{Frame~6: Fermionic observable}}

Consider a simply connected domain $\Omega\subset\mathbb C$ and for $a>0$, let $\mathbf \Omega$ be the largest connected component of $a\mathbb Z^2\cap \Omega$. Consider  $n$ vertices $x_1,\dots,x_n$ of $\mathbf \Omega$, and $n$ faces $f_1,\dots,f_n$ of $\mathbf \Omega$ such that $f_i$ is bordered by $x_i$ for every $1\le i\le n$. Choose $n$ disjoint {\em cuts} $\ell_1,\dots,\ell_n$, i.e.~families of dual edges $(e_i^*(j))$ forming self-avoiding paths in the dual from the unbounded face to the center of $f_i$. Define the {\em disorder operator} $\mu_\ell$ for a cut $\ell$ as the observable that effectively switches the coupling constants of the edges $e_i(j)$ associated with the $e_i^*(j)$ in the cut (it can be written as a product of terms of the form $\exp[-2\beta\sigma_x\sigma_y]$ over edges appearing in the family of edges $\{e_i(j):i,j\}$). Then, the {\em order-disorder} correlations are given by the formula
\begin{equation}
F_\Omega^{(a)}(x_1,f_1,\dots, x_n,f_n):=\langle\sigma_{x_1}\mu_{\ell_1}\dots\sigma_{x_n}\mu_{\ell_n}\rangle_{\mathbf \Omega}.
\end{equation}
Let us mention that these quantities can be expressed in terms of correlations of Grassmann variables in the Schultz-Mattis-Lieb representation \cite{SchMatLie64}. 

Smirnov introduced a {\em fermionic observable} $f_\Omega^{(a)}$ defined at centers of edges $\{x,y\}$ of $\mathbf{\Omega}$ that can be written as a linear combination (with complex coefficients) of the $F_\Omega^{(a)}$ with $x_1$ equal to $x$ or $y$, and $f_1$ to one of the two faces bordered by $\{x,y\}$. The details of the definition are unimportant here and the take-home message is that Chelkak and Smirnov proved that the limit (as $a$ tends to 0) of these fermionic observables is conformally covariant. 
\end{framed}

The conformal covariance of the fermionic observable should be understood as the first brick among the conformal covariance results of spin-spin, energy-energy correlations, and even of the conformal invariance of interfaces. Let us mention that these results require substantial additional ideas compared to \cite{Smi10,CheSmi12}. In fact, conformal covariance/invariance of virtually all quantities one may be interested in the 2D Ising model can be recovered today.

The proof of the theorem relies on the observation that $f_\Omega^{(a)}$ is the solution of a {\em discrete} version of a Riemann-Hilbert boundary value problem. More precisely, the function can be proved, via combinatorial arguments involving the van der Waerden high-temperature expansion, to be {\em preholomorphic} (see Frame~7), and to satisfy certain boundary conditions. These special features are connected to the integrability of the model. From general principles on preholomorphic functions, the limit as $a$ tends to 0 of these objects must be the holomorphic solution of a continuum Riemann-Hilbert boundary value problem, which can be computed and proved to be conformally covariant. Such reasoning has been used in several existing proofs of conformal invariance, for instance for dimers or Bernoulli site percolation on the triangular lattice. It has created an explosion of results in the field as many quantities can be proved to converge using a similar strategy.

\begin{framed}
\centerline{\textbf{Frame~7: Preholomorphic observables}}

The notion of preholomorphic function on a planar graph $G$ appeared implicitly in the work of Kirchhoff on electrical networks \cite{Kir47}. It was explicitly linked to holomorphicity in the work of Isaacs \cite{Isa41,Isa52}, in which the author proposed to discretise the Cauchy-Riemann equation to get to the definition (on the square lattice)
\begin{equation}
F(NW)-F(SE)=i[F(NE)-F(SW)],
\end{equation}
where $NW$, $SW$, $SE$, and $NE$ are the four corners found in counterclockwise order around each face, when starting from the top left vertex.

The properties of preholomorphic functions have been the object of a renewed interest with the emergence of the question of conformal invariance in connections to boundary value problems. Indeed, general theorems stating that preholomorphic functions satisfying certain boundary value conditions converge when taking finer and finer meshsize to holomorphic solutions of the continuum version of the boundary value problem took a central place in the theory.

In the case of the Ising model, the complexity of the boundary value problem (involving a condition on the argument of the fermionic observable) pushed Smirnov to introduce a stronger notion of preholomorphicity, called $s$-{\em holomorphicity}, which is also satisfied by fermionic observables. The advantage of this notion is that it enables one to define the imaginary part of the primitive of the square of the observable, which roughly speaking becomes the discrete solution of a Dirichlet boundary value problem, a much more tractable problem for which convergence (when $a$ tends to 0) can be proved very elegantly.
\end{framed}

\subsection{Towards universality of the 2D Ising model}\label{sec:8.3}

As mentioned in Section~\ref{sec:4.2.2}, the large-scale properties of the critical Ising model should not depend on the precise properties of the underlying graph. With the tremendous successes that have been achieved over the years in the case of the Ising model on $\mathbb Z^2$ and more generally on planar graphs, it is natural to test the validity of the universality hypothesis in this context. Several advances have been made in this direction in the last fifteen years.

The first impressive progress can be found in the work of Chelkak and Smirnov themselves \cite{CheSmi12}. They observed that the preholomorphicity argument leading to conformal invariance can be articulated naturally in the setting of so-called {\em isoradial graphs}. An isoradial graph is an {\em embedding} of a graph $G$ in the plane such that every face of the graph is inscribed in a circle of radius 1.
In this context, one may define special coupling constants $J_{x,y}$ depending on the graph in such a way that $\beta_c=1$ and that the fermionic observable is naturally preholomorphic on this graph. Then, the strategy of Chelkak and Smirnov on the square lattice applies to isoradial graphs with the same conclusions. Note that this result can be understood as a universality result on the graph (isoradial graphs are a fairly large family of planar graphs, even though not fully general), but that the choice of $J_{x,y}$ is {\em determined} by the embedded graph itself. Moreover, a striking feature of this theorem is that no transitivity or quasi-transitivity is required for this to work. 

In recent developments, Chelkak generalised the conformal invariance result to a wider class of Ising models, namely those defined on planar locally-finite doubly periodic weighted graphs $(G ,J)$,
i.e.~weighted graphs which are invariant under the action of some lattice $\Lambda\approx \mathbb Z\oplus\mathbb Z$ (in such case $G/\Lambda$ is a finite graph embedded in the torus).  For such models, Chelkak proved in \cite{Che20} that there exists  an embedding in the plane, called an {\em $s$-embedding}, with the property that the scaling limit of the critical model defined on this embedding is conformally invariant.

This result is a strong indication of universality for planar graphs. Now what happens beyond planar graphs?
The universality conjecture asserts that the scaling limit depends on the {\em large scale} geometry of the graph (for instance a planar Euclidean geometry). In particular, one may consider the graph obtained with the vertex-set $\mathbb Z^2$ and edge-set given by pairs of vertices at a distance at most $R$ of each other. This model, called the {\em finite-range} model on $\mathbb Z^2$, should have a behaviour that is similar to the nearest-neighbour case as it is ``almost planar''. 
The additional difficulty is that non-planarity immediately breaks the integrability of the system.
The universality of such Ising models has been investigated in two different directions. 

First, one may consider finite-range models that are perturbations of the nearest-neighbour integrable case, meaning that non-nearest neighbour interactions are very weak, i.e.~that $J_{x,y}$ is small when $1<\|x-y\|_2\le R$. Using the Schultz-Mattis-Lieb Grassmann representation \cite{SchMatLie64} of the nearest neighbour case, one may express the partition function and more generally the energy-energy and spin-spin correlations in terms of Grassmann variables, and therefore at the end in terms of the nearest-neighbour model. Using an elaborate multi-scale analysis and studying the renormalisation of parameters induced by this multi-scale analysis, Giuliani-Greenblatt-Mastropietro derived in \cite{GiuGreMas12} the large-scale behaviour of energy-energy correlations in the full plane. While the previously mentioned renormalisation schemes in dimension 4 were enabled by the fact that the model is a small perturbation of the discrete GFF (which is a gaussian process), the two-dimensional case relies on a similar connection, this time to the n.n.f.~Ising model on $\mathbb Z^2$ (which has a Grassmannian structure). As a consequence, the strategy suffers from the same limitations as the 4D case in the sense that it is restricted to small perturbations of the n.n.f. Ising model on $\mathbb Z^2$.

A totally different approach explaining the emergence of planarity in finite range Ising models was proposed in \cite{AizDumTasWar19}
 based on the random current representation. The underlying idea relies on the fact that thanks to the switching lemma, intersection properties of random currents with sources are related to the structure of $n$-point correlations in the model. Yet, the intersection properties of long paths on the graph induced by $\mathbb Z^2$ and the edges between vertices at a distance $R$ of each other resemble the ones that can be obtained for planar graphs. As an example of a possible application, one can obtain that spin-spin correlations on the boundary of a domain $\Omega$ have a Pfaffian structure, a result which is specific to the universality class of the 2D Ising model. More precisely, for any collection of points  $x_1=(k_1,0),\dots,x_{2n}=(k_{2n},0)$  satisfying $k_1<k_2<\dots<k_{2n}$ on the boundary of the upper half-plane $\mathbb H:=\mathbb Z\times\mathbb Z_+$, 
 \begin{equation} \label{Pf_cor}
 \langle \sigma_{x_1}\cdots\sigma_{x_{2n}}\rangle_{\mathbb H,\beta_c} \ = \
\mathrm{Pfaff}_n \big( \big[\langle \sigma_{x_i}\sigma_{x_j}\rangle_{\mathbb H,\beta_c}\big]_{1\le i<j\le 2n} \big)\big[1 +o(1)\big]
\, ,
\end{equation}
where $o(1)$ is a function of the points $x_1,\dots,x_{2n}$ which tends to zero for configuration sequences with  $\min\{|x_i-x_j|:1\le i<j\le 2n\}$ tending to infinity.

This is, to the author's knowledge, the first property witnessing the 2D Ising universality class that can be obtained in a level of generality that is not restricted to planar graphs and their perturbations. Also, the proof relies on the key properties of the Ising model that one would like to use: the $\pm$ spin symmetry (entering the story through the use of the random current representation) and the large scale planarity of the underlying graph (which for finite range models on $\mathbb Z^2$ is the reason behind the ``almost'' intersection properties of long paths).  The trade-off is that full conformal invariance of this family of models is still out of reach.

\subsection{Conformal bootstrap in 3D Ising model}\label{sec:8.4}

At this point, we already mentioned that the 1D Ising model was trivially solved in the original paper of Ising~\cite{Isi25}, and that it took 20 more years to achieve a solution of the 2D Ising model \cite{Ons44}. We also saw that the model in dimensions 4 and higher is much simpler as its large-scale properties should be Gaussian. This singles out 3D as the remaining challenging dimension. To the best of our knowledge \cite{VisPorRapdaL22}, it is not known whether the model is integrable or not. This is particularly problematic as the third dimension is probably the most relevant one physically (for instance the model should be in the universality class of liquid--vapour systems, and totally anisotropic magnets).

In recent years, a striking progress has been made on the physics side using the so-called {\em conformal bootstrap}. 
A conformal field theory (CFT) is characterised by the correlation functions $\langle --\rangle$ of an infinite number of {\em local operators} $\mathcal A(x)$, which in the case of Ising should be understood as the objects obtained by taking the limit of random variables defined in terms of spins next to a given position of space. For example, the scaling limit of spin and energy observables $\sigma_x$ and $\varepsilon_{\{x,y\}}=\sigma_x\sigma_y-\langle\sigma_x\sigma_y\rangle$ give such local operators in the case of Ising, but one may think of more complicated ones, such as the scaling limit of (products of) the gradient $\sigma_{x+y}-\sigma_x$ of the spins. 

Conformal invariance already forces huge constraints on the correlations of operators in the theory. Oversimplifying slightly, for scalar local operators there must exist exponents $\bm\Delta_{\mathcal A}$ and coefficients $f_{\mathcal A\mathcal B \mathcal C}$ such that  
\begin{align}
\label{eq:2pt}\langle\mathcal A(x)\mathcal A(y)\rangle&=\frac1{\|x-y\|_2^{\bm\Delta_{\mathcal A}}},\\ 
\langle\mathcal A(x)\mathcal B(y)\mathcal C(z)\rangle&=\frac{f_{\mathcal A\mathcal B\mathcal C}}{\|x-y\|_2^{\bm\Delta_\mathcal A+\bm\Delta_\mathcal B-\bm\Delta_\mathcal C}\|y-z\|_2^{\bm\Delta_\mathcal B+\bm\Delta_\mathcal C-\bm\Delta_\mathcal A}\|z-x\|_2^{\bm\Delta_\mathcal C+\bm\Delta_\mathcal A-\bm\Delta_\mathcal B}}
\end{align}
(in \eqref{eq:2pt}, we adopted without loss of generality the normalization of $\mathcal A$ that makes the constant in the numerator equal to 1).
The exponents and coefficients depend a priori on the CFT, but a striking feature is that there exists a way, called the {\em conformal block decomposition}, to express multi-point correlations of local operators in terms of three-point functions by gluing points together using the so-called operator product expansion. This theoretically shows that all the information in a CFT can be encoded in terms of the $\Delta_\mathcal A$ and the $f_{\mathcal A\mathcal B\mathcal C}$. Of course, determining these coefficients is very difficult.

While in 2D this was done in the eighties, the analogous question remains widely open in 3D. Nevertheless, one can proceed in a slightly different way by asking which choices of these quantities can lead to a consistent CFT. This approach, called the {\em conformal bootstrap}, was shown to be amazingly powerful in 3D. The underlying idea is that one is facing an infinite family of consistency relations coming from different ways of applying the conformal block decomposition (which is not unique). For instance, one may start with $\langle\mathcal A(x_1)\mathcal A(x_2)\mathcal A(x_3)\mathcal A(x_4)\rangle$ and proceed by gluing first $x_1$ and $x_2$ or, on the contrary, $x_3$ and $x_4$. This leads to two decompositions of the same object as a linear combination (with positive coefficients in the Ising case) of known objects called the {\em conformal blocks}. Equalling these two decompositions, one ends up with constraints on the possible exponents. 

There is a priori no reason to be able to determine the critical exponents as the unique values satisfying a (finite) number of constraints thus obtained. Indeed, the set of possible values may not shrink when considering more and more conditions, but it happens that in the case of the Ising model, the region of the plane for possible critical exponents $(\bm\Delta_\sigma,\bm\Delta_\varepsilon)$ for the spin and energy local operators
can be reduced drastically, to a point where estimates -- namely  $(\bm\Delta_\sigma,\bm\Delta_\varepsilon)=(0.5181489(10), 1.412625(10))$ -- using this bootstrap technique become way better than Monte-Carlo simulations. We refer to \cite{RatRycTonVic08,KosPolSim14,ElSPauPolRycSim12} for some of the original papers and \cite{Sim17} for a review of the most recent progress in this very exciting area of modern theoretical physics.

Let us conclude that even if one may use conformal bootstrap to {\em exactly} identify the critical exponents, this would leave the question of proving that the critical 3D Ising model {\em indeed} converges to a CFT widely open. In some sense, getting sufficient information on the possible scaling limits and proving that these scaling limits indeed exist are two almost entirely disjoint questions even though, of course, one may hope that information on the former question would help answer the latter.

\section{A tail to this story}\label{sec:9}

The Ising model has always played the role of a locomotive in the developments of statistical physics. Its central place and incredible properties turn it into an amazing playground for both mathematicians and physicists. As a consequence, during most of its history novel techniques were developed to solve problems on it, which later led to whole independent fields of mathematical physics (integrable systems, graphical representations, rigorous renormalisation methods, etc). 

Let us mention several long-standing problems remaining widely open for this model. At the top of the list, universality of the 2D behaviour (see Section~\ref{sec:8.3}), critical properties of the  3D model (see Section~\ref{sec:8.4}), and the roughening phase transition (see Section~\ref{sec:5.3}) are among the most important unsolved puzzles. Solving them will probably require the development of new techniques that will again, through cross-fertilisation, benefit the whole field of statistical mechanics.

\begin{ack}
We thank all  our coauthors for the wonderful years of joint research, both past and future. We also wish to thank D.~Cimasoni, T.~Gunaratnam, D.~Krachun, I.~Manolescu, R.~Panis, V.~Tassion, and Y. Velenik
 who, without knowing its ultimate aim, took the time to take a look at this review and to give the author feedback. We thank T.~Hutchcroft for bringing \cite{Jon} to our attention.
\end{ack}

\begin{funding}
This work has received funding from the European Research Council (ERC) under the European Union's Horizon 2020 research and innovation programme (grant agreements No.~757296). The author acknowledges funding from the NCCR SwissMap, the Swiss FNS, and the Simons collaboration on localization of waves.
\end{funding}


\end{document}